\documentclass[12pt,twosided]{amsart}  %leqno =  numeracion de las formulas a la izquierda no a la derecha 
                                            
\usepackage{amsmath,amscd,amssymb,amsfonts, mathrsfs}
\usepackage[all]{xy}
\usepackage[usenames,dvipsnames]{xcolor}

\newcommand{\sF}{{\mathcal F}}

\newcommand{\sH}{{\mathcal H}}

\newcommand{\sM}{{\mathcal M}}

\newcommand{\sO}{{\mathcal O}}

\newcommand{\sW}{{\mathcal W}}
\newcommand{\sX}{{\mathcal X}}

\newcommand{\Hs}{{\mathscr H}}
\newcommand{\C}{{\mathbb C}}

\renewcommand{\H}{{\mathbb H}}

\newcommand{\M}{{\mathbb M}}

\newcommand{\PP}{{\mathbb P}}
\newcommand{\Q}{{\mathbb Q}}

\newcommand{\Z}{{\mathbb Z}}
\newcommand{\Dif}{\Omega}
\newcommand{\dif}{\omega }

\newcommand{\lra}{\longrightarrow }

\newcommand{\Gr}{{\rm Gr}}

\newcommand\codim{\operatorname{codim}}
\newcommand\Ker{\operatorname{Ker}}
\newcommand{\rank}{\operatorname{rank}}
\newcommand{\im}{\operatorname{Im}}
\newcommand{\Cok}{\operatorname{Coker}}

\newcommand{\Sing}{\operatorname{Sing}}
 
\newcommand{\depth}{\operatorname{depth}}
\newtheorem{thm}{Theorem}
\newtheorem{lem}[thm]{Lemma}
\newtheorem{cor}[thm]{Corollary}
\newtheorem{prop}[thm]{Proposition}
\newtheorem{remark}[thm]{Remark}

\newtheorem{defn}[thm]{Definition}
\newtheorem{ej}[thm]{Example}

\def\pf{\noindent{\it Proof. }}
\def\qed{\rm Q.E.D.}

\begin{document}

%\title[\resizebox{4.0 in}{!}{Variation of Mixed Hodge Structures associated to an equisingular \\ one-dimensional family \\ of Calabi-Yau 3-folds}]
%{Variation of Mixed Hodge Structures associated to an equisingular \\ one-dimensional family \\ of Calabi-Yau 3-folds} 

\title[VMHS of an equisingular one-dimensional family of CY threefolds]{Variation of  Mixed Hodge Structures  associated 
to an equisingular 
one-dimensional family \\
of Calabi-Yau threefolds}
 %\title{for}
\thanks{ Mathematics Subject Classification (MSC2010) : 14D07 ~(primary),~14J32,~14C30,~14Q15.}
\author[Nieto, I  and Del Angel, P.L.]{Isidro Nieto-Ba\~nos ~~~~\and  Pedro Luis del Angel-Rodriguez}  
\thanks{ 
The first author acknowledges partial support from  CIMAT   as well as from the CMO congress `` Primer Congreso Nacional de Geometr{\'i}a Algebraica''
in the year 2016 in Oaxaca where part of these 
results were presented. Both  authors had fruitful discussions at various stages of this work with  D. Van Straten, Ch. Peters, J. Carlson, X. Gomez-Mont and H. Kanarek and thank them for it. We acknowledge  partial support from  CONACyT Grant 0181730} 
\email{\tt~nietoisidrorafael@yahoo.com~, luis@cimat.mx } 
\address{ ~C.P.~36023,~Guanajuato,~Gto.,~M\'exico and  CIMAT, A.C., Jalisco, S/N, Guanajuato,~Gto.,~M\'exico.}

\begin{abstract}

We study the variations of mixed Hodge structures 
(VMHS) associated to a pencil~ $\sX$ of
equisingular hypersurfaces of degree $d$ in $\PP^{4}$ with only ordinary double points as singularities, as well as 
the variations of Hodge structures (VHS) associated to the desingularization of this family  $ \widetilde{\sX}$.
The notion of a set of singular points being in {\it homologically good position} is introduced and, 
by  requiring that  the subset of nodes in (algebraic) general position is also in homologically good position, 
we can extend  Griffiths description of the $F^2 $-term  of the Hodge filtration of the desingularization to this case, where we 
can also determine the possible limiting mixed Hodge structures (LMHS).
~The particular pencil $\sX$ of  quintic hypersurfaces with $100$ 
singular double points with $86$ of them in (algebraic) general position 
which served as the starting point for this paper is treated with particular attention.
%including a few remarks concerning its Picard-Fuchs operator.

\end{abstract}                           %= pie de notas 

\maketitle
\date{}

\section*{Introduction}
In 1941 W.V.D. Hodge proved that the complex de Rham cohomology $H^k(X,\C)$ of every compact K\"ahler manifold,  splits as a
direct sum of spaces $H^{p.q}(\cong H^q(X, \Dif_X^p))$, where $p+q=k$, called 
nowadays the Hodge decompositon of $H^k(X,\C)$ (see \cite{Hodge1}). The pair $(H^k(X,\Z), \{H^{p,q}\})$ is
called a (pure) Hodge structure of weight $k$. All varieties will be considered algebraic and defined
over the complex numbers $\C$. Unless otherwise stated we will be consistent with Deligne's notation (see \cite{deligne2}).
\\
 \\
Another way of looking at a Hodge structure is to consider the associated Hodge filtration 
$
{\displaystyle F^jH^k(X,\C) \stackrel{def}{=}\oplus_{p\ge j}H^{p,q}}
$
and the pair $(H^k(X,\Z), \{F^jH^k(X,\C)\})$. \\ \\
If $X\subset\PP^{n+1}$ is a hypersurface, then the only interesting cohomology group is 
$H^n(X,\C)$ and because of Lefschetz' theorem, we only need to consider the so called primitive cohomology $PH^n(X,\C) = \{\eta\in H^n(X,\C)\; |\; \eta\cdot H=0\}$, where
$H$ is the class of a hyperplane section on the corresponding projective space.
\\ \\ 
Griffiths studied the (pure) Hodge structure of smooth projective hypersurfaces $X$  and gave a description of it in terms of 
its Jacobian ring (see \cite{griffiths1}). More precisely, let $X=V(f)\subset\PP^{n+1}$ be a smooth hypersurface of degree $d$ and let
\begin{eqnarray}
 \label{smooth}
\Hs_k(X) \stackrel{def}= \left\{\left[\frac{P\Dif}{f^k}\right]  \; \in A^{n+1}_k\; \mbox{mod} \; dA^{n}_{k-1}\; |\; \mbox{deg}(P)=kd-(n+2)\right\},
\end{eqnarray}
where $A_k^j$ denotes the space of rational $j$-forms on $\PP^{n+1}$ with a pole of order $k$ along $X$ and
$ \Dif= {\displaystyle \sum_{i=0}^{n+1} (-1)^i x_i dx_0\wedge\ldots \wedge \widehat{dx}_i\wedge\ldots\wedge dx_{n+1}} $. Then we have a commutative diagram
\begin{eqnarray}
 \label{GHodge-F}
 \xymatrix{
0 \ar[d] \ar@{^{(}->}[r] & \Hs_1(X)\ar[d]_{\cong} \ar@{^{(}->}[r] & \cdots \ar@{^{(}->}[r] & \Hs_n (X)\ar[d]_{\cong}\ar@{^{(}->}[r] & \Hs_{n+1}(X)\ar[d]^{\cong} \cr
0 \ar@{^{(}->}[r] & F^nPH^n(X,\C) \ar@{^{(}->}[r] & \cdots \ar@{^{(}->}[r] & F^1PH^n(X,\C) \ar@{^{(}->}[r] & F^0PH^n(X,\C),
}
\end{eqnarray}  
\ \\ 
where the horizontal arrows in the first line correspond to the natural inclusion given by multiplication by $f$. Moreover, if
$J(f)$ is the Jacobian ideal of $f$ and $R_f \stackrel{def}{=}\C[X_0,\ldots,X_n]/J(f)$ is the Jacobian ring of $f$, then the above identification induces isomorphisms between $(R_f)_{(k+1)d-n-2}$ and $PH^{n-k,k}(X,\C)$. 
\\ \\
For singular varieties, Deligne developed in 1971 the theory of mixed Hodge structures (see \cite{deligne2}), which involves in general the existence of a good desingularization due to Hironaka. 
\\ \\
Griffiths and others have tried to give an alternative description for the mixed Hodge structure of a singular variety in some cases.
The most important case for us is that of a singular projective hypersurface on the projective space with isolated singularities, 
the simplest of which is only nodes as singularities. For hypersurfaces of dimension less than or equal to 3,
Griffiths  \cite{griffiths2} (who considers three dimensional hypersurfaces with one ordinary double point) and later on Steenbrik \cite{steenbrik2} (who considers surfaces with isolated singularities) gave 
a description of the relevant cohomology group of its proper transform under normalization in terms 
of the Jacobian ring of the polynomial defining it. More precisely, let $X=V(f)\subset\PP^{n+1}$ be a hypersurface of degree $d$, 
and assume its singular locus $\Sigma$  consist of ordinary double points. Let $\widetilde{X}$ be its proper transform under normalization. 
 \label{fadjoint}
 If we define the $\Hs_k$ as before as well as the vector space
\begin{eqnarray*}
\Hs_2^1(X)  \stackrel{def} =   \left\{\left[\frac{P\Dif}{f^2}\right]  \; \in A^{n+1}_2\; \mbox{mod} \; dA^{n}_{1}\; |\; \mbox{deg}(P)=2d-(n+2)\; \mbox{and} \; P(Q)=0 \; \forall Q\in\Sigma \right\}, 
\end{eqnarray*}
 given by the {\it first adjunction condition on}  $ A^{n+1}_2$, then if $|\Sigma|=1$, we get a partial generalization of commutative diagramm (\ref{GHodge-F}) for the $n-$th cohomology of the primitive part of $\widetilde{X}$, namely
 \begin{eqnarray}
 \label{FA-condition}
 \xymatrix{
0 \ar[d] \ar@{^{(}->}[r] & \Hs_1(X)\ar[d]_{\cong} \ar@{^{(}->}[r] & \Hs_2^1(X) \ar[d]_{\cong}  \cr
0 \ar@{^{(}->}[r] & F^nPH^n(\widetilde{X},\C) \ar@{^{(}->}[r] & F^{n-1}PH^n(\widetilde{X},\C).
}
\end{eqnarray}
 Where the horizontal map in the arrow below corresponds to the natural map
\begin{eqnarray*}
 \Hs_1(X) \longrightarrow  \Hs_2^1(X)
\end{eqnarray*}

\noindent
given by ${\displaystyle \frac{P\Dif}{f} \mapsto \frac{f P\Dif}{f^2}}$. A direct generalisation of (\ref{FA-condition}) for $|\Sigma| \geq 2$ is not so straightforward as it may seem
and we show  in  theorem \ref{taut} part  2) that to assume that all the points of $\Sigma$  are in (algebraic) general position is
not enough, so further imposing the condition that $\Sigma$ is a hg set (see  our definition \ref{hg})  gives the  expected generalization as we have proven  in  corollary \ref{precise} and in corollary \ref{cor:rank}.
\\ \\
If we now consider a smooth family $\pi: \sX\lra B\subset \PP^1$, over a Zariski open set $B$, then on every fiber $X_t$ one has a Hodge structure
$(H^n(X_t, \Z),\{F^pH^n(X_t,\C)\})$
and the Hodge filtration extend to a global filtration $\sF^p\sH^n$, where 
$\sH^n \stackrel{def} = R^n\pi_*\C\otimes\sO_B$.
It is well known that the monodromy of the family gives rise to a connection, called the Gauss-Manin connection (or shortly GM)
\begin{eqnarray*}
\nabla: \sH^n\lra \sH^n\otimes\Dif_B
\end{eqnarray*}
which is compatible with the Hodge filtration. More explicity, the GM-connection satisfies 
the Griffiths transversality condition (called also  the horizontality condition)
\begin{eqnarray*}
\nabla: \sF^p\sH^n\lra \sF^{p-1}\sH^n\otimes\Dif_B.
\end{eqnarray*}

Recall that any polarized VHS $\H$ of weight $k$ on $B$ induces a map from $B$ to the classifying space $D$ of polarized Hodge structures of weight $k$, which can be seen as a Zariski open set on a projective variety $\check{D}$ parametrizing flags $F^{k}\subset\cdots\subset F^1\subset V_{\C}$ of nationality $(f^k,\cdots, f^1)$ satisfying the first Riemann-Hodge bilinear relation, where $f^j\stackrel {\text{def}} = \dim F^j$. Since $\check{D}$ is projective, the map
$B\to D$ induced by $\H$ can be extended to a map $\PP^1\to \check{D}$, in particular to any point $p\in\PP^1\backslash B$ we can associate a filtration
$F^{k}(p)\subset\cdots\subset F^1(p)\subset V_{\C}$ satisfying the first Riemann-Hodge bilinear relation. There is no reason for this filtration to satisfy the second Riemann-Hodge bilinear relation and in general it will not, so   a priori there  does not exist a  polarized Hodge structure of weight $k$ at $p$. However, in a small analytic neighborhood $U$ of $p\in\PP^1\backslash B$, the local monodromy $\pi_1(U,b)\cong\Z$ regardless of the choice of $b\in U$, and the generator of this group induces a linear transformation $T$ on $\H_b$, compatible with the Hodge filtration, called the \emph{monodromy operator}. $T$ can be thought of as a linear transformation of $V_{\C}$ and as such, it will be compatible with the filtration $F^{k}(p)\subset\cdots\subset F^1(p)\subset V_{\C}$. If we write $T=T_s \circ T_u$, where $T_s$ and $T_u$ are the semisimple  and the unipotent part of $T$ respectively. It can be proved  (see \cite{peters} Chapter 11 theorem 11.8 (monodromy theorem) and  lemma-definition 11.9)   or \cite{Kulikov} - monodromy theorem \S 9.1 ) that $T_s^m=1$ for some integer $m$ and the least $l$ such that  $(T_u - \text{Id})^{l}=0$ is less than or equal to $m+1$. The triple 
$(W_{\bullet}, F^{\bullet}, V)$ (or simply $(W_{\bullet}, F^{\bullet})$ whenever $V$ is clear from the context), defines a MHS on $V$, called the \emph{limit MHS} at $p$ ( in the sense of Schmid), where $W_{\bullet}$ is the monodromy weight filtration on $V$ associated to $N=\text{log}(T_u)$ (see \cite{griffiths4} page 255, \cite{griffiths3} pages 106-107 and  lemma \ref{lem:monodromyweight}). Finally, since the GM-connection associated to the VHS $\H$ satisfies Griffiths transversality, the nilpotent operator $N$ induces linear maps $N_j:H^{j,k-j}\to H^{j-1,k-j+1}$, where $H^{j,k-j}\stackrel{\text{def}} = F^j(p)/F^{j+1}(p)$. The previous discussion  is true for the case of $PR^k\pi_*\C$,  the primitive part of the $k$-th higher direct image  of $\C$.
 \\ \\
 From now on, unless explicitly stated otherwise,  $X\subset\PP^{n+1}$ will be a singular hypersurface with singular locus $\Sigma$, $\pi:\widehat{\PP}^{n+1}\to \PP^{n+1}$ will denote the blow up of $\PP^{n+1}$ along $\Sigma$ and $\widehat{\Sigma}$ will be the exceptional divisor on $\widehat{\PP}^{n+1}$, moreover $\widetilde{X}$ will be the strict transform of $X$ and $\widetilde{\Sigma}=\widehat{\Sigma}\cap \widetilde{X}$. 
 \\ \\
The paper is organized as follows. In section 1, we  generalize  the classical definition of adjointness of $\Hs_2^1$ to  isolated singularities of higher order, denoted as s-adjointness in definition \ref{sadj}. The most important result in this section is given by  proposition \ref{OmegaX}, which is a sheaf theoretic formulation  of the notion of s-adjointness with pole order conditions and a partial description of $\Hs^k(\widetilde{X})$ in terms of rational forms on $\PP^4$, at least when the singular locus of $X$ consist of ordinary double points. In section 2, we define the notion of homologically good sets and study its relations to the notion of points in algebraic general position. The central results are theorem \ref{taut}, 
 generalizing  diagramm (\ref{FA-condition}) above,  together with corollaries  \ref{precise} and \ref{cor:rank}. In section 3, recalling  the 
 definition of generalized Hodge numbers, the main result is  given by proposition \ref{ghn}  which computes the generalized Euler characteristic polynomial of
 $X$ and $\widetilde{X}$ using the techniques  introduced in \cite{danilov}. In section 4, using the technique of cubical hyperresolutions of \cite{Navarro} the main results are: First, the computation of the Mixed Hodge structure of a nodal threefold (proposition \ref{GrF}) applied in
 example  \ref{ej:quintic},
  where we actually compute the number of points in algebraic general position. The other important result in  this section is given in proposition \ref{defail} and remark \ref{defail2} by finding an exact relation between the defect and the failure of $\Sigma$ to impose linearly independent conditions on polynomials of degree $2d-5$. In  section 5, we consider the VHS associated to the Lefschetz pencil of the desingularizations and compute the possible weight filtrations  corresponding to the limit MHS in theorem \ref{lmhs}. Another important and natural result is that  the
 VMHS associated to an equisingular pencil of nodal threefolds  is  indeed a geometric and admissible VMHS in the sense of \cite{stzu}, stated and proved in proposition \ref{locsys}. 
 
\section{Generalized adjointness conditions}

Let $ \Omega_{\PP}^4 (kX) $ be the sheaf on $ \PP^4 $ of   four-rational forms  with a pole of order
$k$ along  the hypersurface $ X$ or shortly $ \Omega_{\PP}^4(k) $. Then it follows that 
$H^0(\PP^4, \Omega_{\PP}^4(k)) = A_k^4(X)$. For a polynomial $F$ we denote by $ \mu_p( F) $ 
the multiplicity of  $F$ in $P$ ( see \cite{Ha}). Given a subset $T\subset X$ let us denote by
$ \Omega_{\PP}^4( kX, s T)$ the subsheaf of $\Omega_{\PP}^4(k) $ of four-rational forms with a pole of order 
$k$ on  $X$ and multiplicity  at least  $s$ on every point $P\in T$. 

\begin{defn}
 Given $f \in \C[y_0,\ldots, y_n]$ the $s$-adjoint condition on $f$ relative to $T$ is given by
 $\mu_p(f) \ge s$ for all  $  p\in  T$. Note that $s = 1$  if and only if $ T \subset V(f)$.
\end{defn}

\begin{defn}
\label{sadj}
The space of four-rational forms with poles of order $k$ along $X$ and s-adjoint to $T$ is defined as follows:
\[ A_{k}^4(X, s T) = \{ \psi \in  A_k^4(X) | \psi = \frac{h \Dif}{f^{k}}, \,\,  h \,\,  \mbox{is } s-\mbox{adjoint relative to} \, \, T \}\]
\end{defn}

\noindent In particular, if $T=\Sigma = \mbox{Sing}(X)$, it follows that $ H^0(\PP, \Dif_{\PP}^4(kX, s\Sigma)) ={ A_k^4}(X_t, s\Sigma) $.
Clearly $ s \leq d = \mbox{deg}(f) $. We have already defined the vector space $\Hs_2^1$ following  \cite{griffiths1}'s notation, and it is clear that
$\Hs_2^1=A^4_2(X_t, \Sigma)/dA^4_1(X)$. If $\Sigma$ consists of ordinary double points  
then trivially  $ d A_{1}^3(X) \subset 
A_{2}^4(X, \Sigma)$ but in general it is  not even possible to compare $ d A_{k-1}^3(X)$ with $A_{k}^4(X, (k-1) \Sigma) $. 
Hence   we can define the following quotient:
$ {\Hs}_k^s ={\displaystyle \frac{ A_{k}^4(X,s \Sigma)}{ d A_{k-1}^3 \cap A_k^4(X,s \Sigma)}} $, which is the vector space of top
rational forms with poles of order $k$ along $X$ and satisfying the $s$-adjoint condition relative to $\Sigma$, naturally 
generalising one-adjointness relative to $\Sigma$ given by \cite{griffiths2}.

\begin{remark}
In this sense given $G $, a finite subset of polynomials, one can generalize the adjointness conditon relative to $T$ if for all
$ h \in G$  the $s$-adjoint condition is satisfied on $h$. 
\end{remark}

 Let us return to the sheaf theoretic version of forms with pole order and adjointness conditions: \\ \\
\textbf{Note:}
We will often write  ${\Dif}_{\PP}^4(k,s) $ as short hand notation for ${\Dif}_{{\PP}^4}^4(k X,s \Sigma )$. Analogously, $\Dif_{\widehat{\PP}}^4(k)$ will stand for $\Dif_{\widehat{\PP}^4}^4(k\widetilde{X})$. 
\\ \\

\begin{prop} 
\label{OmegaX}
 With notation as above, if $\Sigma$ consists of ordinary double points, then for $N= 2k-3$ positive and $ s \geq N$ we have  $ \pi^{*}\big( {\Dif}_{\PP}^4(k,s)\big) \subset \Dif_{\hat{\PP}}^4(k)$.
\end{prop}

\label{mult}
\pf
This is a local straightforward computation.  \hfill \qed   
% % % % % % % % % % % % % % % % % % % % % % % % % % % % % % % % % % % % % 28/12/14 % % % % % % % % % % % % % %
%%%%%%%%%%%%%%%%%%%%%%%%%%%%%%%%%%%%%%%%%%%%%%%%%%%%%%Nie-Ang 23/04/15%%%%%%%%%%%%%%%%%%%%%%%%
\section{Elementary results for nodal hypersurfaces on $\PP^4$}
\label{hyp}

Given a projective variety $X$, we will say that a finite set $T\subset X$ is a set of points in \emph{ algebraic general  position} or shortly 
{\it in general position } if they impose $|T|$ conditions on polynomials of degree $d$ passing through all of them, for all $d\ge 1$.

\begin{lem} \label{CohomCompSupp}
For any scheme $Y$ of dimension $n>0$, any locally free sheaf $H$ of finite rank and any non-singular subvariety $Z\in Y$, we have $H^i_{Z}(Y,H)=0$ for all $i<codim_{Y}(Z)$.
\end{lem}
\pf It follows by excision, since for any $\text{P}\in Y$ with smooth closure $Z$, the local ring $\sO_{\text{P}}$ is a regular local ring of $\depth = \codim_{Y}(Z)$. \hfill \qed
\ \\ \ \\
%Assume now that $X$ is a 3-fold on $\PP^4$ of degree $\ge %5$ with $m$ nodes, 
%where  precisely $l\le m$  are in general position. Then 
The central result in this section is given by the following:
\begin{thm}
\label{taut}
\ \\
\renewcommand{\theenumi }{\alph{enumi}}
\begin{enumerate}
\item  
$\xymatrix{\Hs_1\ar[r]^(0.4){\pi^*} & H^{3,0}(\widetilde{X})}$
is an isomorphism and
\item $\xymatrix{\Hs_2^1\ar[r]^(0.37){\pi^*} & F^2H^3(\widetilde{X}, \C)}$ is injective.
\end{enumerate}
\renewcommand{\theenumi }{\arabic{enumi}}
\end{thm}
\pf 
It is well known that  $H^{3,0}(\widetilde{X})\cong H^0(\widehat{\PP}^4, \Dif^4_{\widehat{\PP}^4}(\widetilde{X}))$,  as already shown  in  \cite{griffiths2}, theorem 10.8 and the first assertion is equivalent to 
$\pi^{*}: H^0(\PP^4, \Dif^4_{\PP^4}(X))\to H^0(\widehat{\PP}^4, \Dif^4_{\widehat{\PP}^4}(\widetilde{X}))$
being an isomorphism.
\\ \\
Let $\xymatrix{U=\PP^4\backslash \Sigma\ar@{^(->}[r]^(0.65){j}& \PP^4}$ and 
$\xymatrix{\widehat{U}=\widehat{\PP}^4\backslash \widehat{\Sigma}\ar@{^(->}[r]^(0.65){\hat{\j }} & \widehat{\PP}^4}$. Then $\pi|_{\widehat{U}}: \widehat{U}\to U$ is an isomorphism and in particular we have a commutative diagram
{\small                       
\begin{eqnarray*}
\xymatrix{
0\ar[r] & H^0_{\Sigma}(\PP^4, \Dif^4_{\PP^4}(X))\ar[r] & H^0(\PP^4, \Dif^4_{\PP^4}(X))\ar[r]^{j^{*}}\ar[d]^{\pi^{*}} & H^0(U, \Dif^4_{\PP^4}(X)) \ar[r]\ar[d]^{\pi^{*}|_U} & H^1_{\Sigma}(\PP^4, \Dif^4_{\PP^4}(X))\cr
0\ar[r] & H^0_{\widehat{\Sigma}}(\widehat{\PP}^4, \Dif^4_{\widehat{\PP^4}}(\widetilde{X}))\ar[r] & H^0(\widehat{\PP}^4, \Dif^4_{\widehat{\PP}^4}(\widetilde{X}))\ar[r]^{\hat{{\j }}^{-1}} & H^0(\widehat{U}, \Dif^4_{\widehat{\PP}^4}(\widetilde{X}))\ar[r] & H^1_{\widehat{\Sigma}}(\widehat{\PP}^4, \Dif^4_{\widehat{\PP}^4}(\widetilde{X}))
}
\end{eqnarray*}
}
where the map ${\pi}^{*}|_U$ is an isomorphism and the rows are exact. 
\\ \\
Since $\widehat{\PP}^4$ is regular in codimension 1 and $\widehat{\Sigma}$ is a disjoint union of exceptional divisors $E_P$ 
above the points $P\in\Sigma$, then
$ H^0_{\widehat{\Sigma}}(\widehat{\PP}^4, \Dif^4_{\widehat{\PP}^4}(\widetilde{X}))
\cong \oplus_{P\in\Sigma} H^0_{E_P}(\widehat{\PP}^4, \Dif^4_{\widehat{\PP}^4}(\widetilde{X}))=0$. 
\ \\ \ \\
On the other hand, since $\Dif^4_{\PP^4}$ is locally free and every point in $\Sigma$ is a regular point in $\PP^4$, by  lemma(\ref{CohomCompSupp}) one also has that $H^1_{\Sigma}(\PP^4, \Dif^4_{\PP^4}(X))=H^0_{\Sigma}(\PP^4, \Dif^4_{\PP^4}(X)) = 0$, so the diagram above becomes
\[
\xymatrix{
0\ar[r] & H^0(\PP^4, \Dif^4_{\PP^4}(X))\ar[r]^{j^{*}}\ar[d]^{\pi^{*}} & H^0(U, \Dif^4_{\PP^4}(X)) \ar[r]\ar[d]^{\pi^{*}|_U} & 0\cr
0\ar[r] & H^0(\widehat{\PP}^4, \Dif^4_{\widehat{\PP}^4}(\widetilde{X}))\ar[r]^{\hat{{\j }}^{*}} & H^0(\widehat{U}, \Dif^4_{\widehat{\PP}^4}(\widetilde{X}))\ar[r] & H^1_{\widehat{\Sigma}}(\widehat{\PP}^4, \Dif^4_{\widehat{\PP}^4}(\widetilde{X}))
}
\]
Since $\pi^{*} \circ j^{*}=\hat{\j}^{*}\circ\pi^{*}$ is an isomorphism, then $\pi^{*}$ is injective and $\hat{\j}^{*}$ is surjective. But $\hat{\j}^{*}$ is injective, therefore $\hat{\j}^{*}$ and $\pi^{*}$ are isomorphisms.
\ \\ \\
As for 2): if we denote by $d $ the total differential  then we have the following  commutative diagram:

\begin{eqnarray}
\label{ctaut}
\xymatrix{ \ar@{->}[d]^{d} A^3_1(X) \ar@{^{(}->}[r]^{{\pi}_1^*} &
 \ar@{->}[d]^{\tilde{d}} A^3_1(\widetilde{X})   \\
          \ar@{->}[d]^{p}   A^4_2(X, \Sigma ) \ar@{->}[r]^{{\pi}^*} & \ar@{->}[d]^{q} A^4_2(\widetilde{X}) \\
             \Hs^1_2(X) \ar@{->}[r]^{\lceil {\pi^\ast} \rceil} & \Hs_2(\widetilde{X}) } 
\end{eqnarray}
where we distinguish the pullback on 3-forms from the pullback on 4-forms through the subindex 1
and
where  $p$ and $q$ are the natural quotient maps. Observe that $\Cok(d) =\Hs^1_2 $  and  similarly $\Cok(\tilde{d}) = \Hs_2(\widetilde{X}) $  and the 
last horizontal arrow
 $\lceil {\pi^\ast} \rceil$ is induced by the universal property of the quotient. 
\\ \\
\noindent {\it Claim 1.} ${\pi}_1^*$ is an isomorphism.

\noindent \pf (of claim 1.) 
\\ \\
\noindent As before, we have a commutative diagram with exact rows:
{\small                       
\[
\xymatrix{
0\ar[r] & H^0_{\Sigma}(\PP^4, \Dif^3_{\PP^4}(X))\ar[r] & H^0(\PP^4, \Dif^3_{\PP^4}(X))\ar[r]^{j^{*}}\ar[d]^{{\pi}_1^{*}} & H^0(U, \Dif^3_{\PP^4}(X)) \ar[r]\ar[d]^{{\pi}_1^{*}|_U} & H^1_{\Sigma}(\PP^4, \Dif^3_{\PP^4}(X))\cr
0\ar[r] & H^0_{\widehat{\Sigma}}(\widehat{\PP}^4, \Dif^3_{\widehat{\PP^4}}(\widetilde{X}))\ar[r] & H^0(\widehat{\PP}^4, \Dif^3_{\widehat{\PP}^4}(\widetilde{X}))\ar[r]^{\hat{{\j }}^{*}} & H^0(\widehat{U}, \Dif^3_{\widehat{\PP}^4}(\widetilde{X}))\ar[r] & H^1_{\widehat{\Sigma}}(\widehat{\PP}^4, \Dif^3_{\widehat{\PP}^4}(\widetilde{X}))
}
\]
}
where the map ${\pi}_1^{*}|_U$ is an isomorphism and, exactly as before, this diagram becomes
\[
\xymatrix{
0\ar[r] & H^0(\PP^4, \Dif^3_{\PP^4}(X))\ar[r]^{j^{*}}\ar[d]^{{\pi}_1^{*}} & H^0(U, \Dif^3_{\PP^4}(X)) \ar[r]\ar[d]^{{\pi}_1^{*}|_U} & 0\cr
0\ar[r] & H^0(\widehat{\PP}^4, \Dif^3_{\widehat{\PP}^4}(\widetilde{X}))\ar[r]^{\hat{{\j }}^{*}} & H^0(\widehat{U}, \Dif^3_{\widehat{\PP}^4}(\widetilde{X}))\ar[r] & H^1_{\widehat{\Sigma}}(\widehat{\PP}^4, \Dif^3_{\widehat{\PP}^4}(\widetilde{X}))
}
\]
and the claim follows.

\noindent {\it Claim 2.} $\pi^{*}$ is injective.

\noindent \pf (Of claim 2.)
\\ \\ Clearly we have a a morphism of sheaves  $\pi^{*}:\Dif^4_{\PP^4}(2,1)\to {\pi}_*\Dif^4_{\widehat{\PP}^4}(2)$
and it  is enough to show injectivity  on the stalk at every point.
\\ \\
If $\widehat{U}$ is an open set on $\widehat{\PP}^4$ whose intersection with $\widehat{\Sigma}$ is empty, then ${\pi}$ is an isomorphism from $\widehat{U}$ to $U={\pi}(\widehat{U})$ and so 
${\pi}^*:\Dif^4_{\PP^4}(2,1)(U)\to {\pi}_*\Dif^4_{\widehat{\PP}^4}(2)(U)=\Dif^4_{\widehat{\PP}^4}(2)(\widehat{U})$ is an isomorphism. Now let us consider an open set $U\subset \PP^4$ containing just the point $P\in\Sigma$, then any rational 4-form $\dif$ on $\PP^4$ with poles of order 2 along $X$ can be written in the form $\frac{F\,d\,z}{(z\cdot z)^2}$, where $z=(z_1, z_2, z_3, z_4)$ are local coordinates and $z\cdot z$ is a local equation defining $X$ on $U$. The form $\dif$ satisfies the first adjoint condition relative to $\Sigma$ in $U$ if and only if $F(P)=0$. But in this case,
\[
{\pi}^*(\dif)(u,v)={\displaystyle \frac{u^3 F(uv,u)\,d\,u\,d\,v}{u^4(1+v\cdot v)^2} } = {\displaystyle \frac{F(uv,u)\,d\,u\,d\,v}{u(1+v\cdot v)^2} }.
\]
Since the zero set of a non constant holomorphic function is a hypersurface on $U$ and ${\pi}$ is a birational morphism, ${\pi}^*(\dif)=0$ if and only if ${\pi}(\widehat{U})\subset V(F)\cap U$ if and only if $F\equiv 0$. 
\ \\ \ \\
\noindent {\it Claim 3}. $\lceil {\pi}^\ast \rceil$ is monomorphism. 

\noindent \pf ( Of claim 3.)
\\ \\
Assume $\overline{\varphi}\in \Hs_2^1$ satisfies $\lceil {\pi}^* \rceil (\overline{\varphi})=0$ and let 
$\varphi\in A^4_2(X, \Sigma )$ be any representative of $\overline{\varphi}$. Then ${\pi}^*(\varphi) = d\, \mathfrak{h}$ for some $\mathfrak{h}\in A^3_1(\widetilde{X})$ and by claim 1 there exist some $\beta\in A^3_1(X) $ such that ${\pi_1^{\ast}}\beta = \frak{h}$, therefore 
\[
{\pi}^*(d\,\beta) = d\, ({\pi}_1^* \beta) = d\, \frak{h} = {\pi}^*(\varphi)
\]
and the injectivity of ${\pi}^*$ implies that $\varphi= d\,\beta\in d\, A_1^3(X)$, i.e. $\overline{\varphi} = 0\in \Hs_2^1$.
\\ \\
The theorem now follows from the fact that
 $ \xymatrix@1{ \Hs_2(\widetilde{X})  \ar@{^{(}->}[r] & F^2H^3(\widetilde{X},\C)} $
as proven e.g. in \cite{griffiths2} proposition 16.3 
 equation  (16.10). \hfill\qed
\\ \\

Let $X\subset \PP^{4}$ be a nodal hypersurface of degree $d$ with $m$ nodes and let $P$ be a node on $X$. Then in an analytic neighborhood $U$ of $P$ in $\PP^4$ we can write $X\cap U=V(z\cdot z)$,
where $ z= (z_1,z_2,z_3,z_4) \in\C^4$ and $z_i = x_i + \sqrt{-1} y_i$  for $  i  \in \{1,2,3,4 \}$. 
With this notation $ z \cdot z = x \cdot x  -y \cdot y +2 \sqrt{-1} ( x \cdot y ) $ and following the technique of continous deformation
 used in~\cite{griffiths2} ( in particular  the notation before equation (15.3)) we can consider
a family $X_{\epsilon}$ of hypersurfaces with $m-1$ nodes that degenerate to $X$, that is to say that
$X-U\cong X_\epsilon - U$.  Note that the three-dimensional real spheres
$ \delta _{\epsilon} = \{x\cdot x=\epsilon, y = 0 \} $ are contained in the $X_{\epsilon} \cap U = \{  z\cdot z = \epsilon \} $ so the 
family of hypersurfaces $ \{ X_{\epsilon } \}$ degenerate to $X $ and the latter is a singular hypersurface with $m$ double points. Observe that there exists a 
$3$-cell $\theta_{\epsilon}(P)$ on $U\cap X_{\epsilon}$ such that $\theta_{\epsilon}(P)\cdot \delta_{\epsilon}(P)=1$
as shown in \cite{griffiths2}. The sphere ${\displaystyle\delta_0(P)\stackrel {def} = \lim_{\epsilon\to 0}\delta_{\epsilon}(P)}$ is contractible to a point in $X$ while the $3$-cell ${\displaystyle \theta_0(P)\stackrel {def} = \lim_{\epsilon\to 0}\theta_{\epsilon}(P)}$ gives a non-zero element of $H_3(X,\Z)\otimes\Q$. However it may happen that $\theta_0(P)$ belongs to the subspace of $H_3(X,\Z)\otimes\Q$ generated by $\{\theta_0(Q)\; |\; Q\text{ is a node of X }, Q\ne P\}$. This motivates the following definition. 
 
\begin{defn}
\label{hg}
Given a nodal hypersurface $X\subset\PP^4$, we will say that a set $T$ of nodes on $X$ is 
\emph{homologically good} (hg) if the 
corresponding set of three-cells $\{\theta_0(P)\}_{P\in T}$, is a $\Q$-linearly independent set of elements in $H_3(X,\Z)\otimes\Q$ and 
$T$ is maximal with this property. 
In particular, there will be a \emph{vanishing cycle} $\delta_0(P)$ for every node in $T$.
\end{defn}

\begin{remark}
Let $X$, $U$ and $X_\epsilon$ be as above, $T$ be a homologically good set of nodes on $X$, $P\in T$ and $\widetilde X$ is the strict transform of $X$ under the blow-up of $X$ on $P$. ~The strict transform of $\theta_0(P)$ does no longer represent an element in the homology of $\widetilde{X}$ since it 
is no longer a cycle in $\widetilde{X}$, and $\delta_0(P)$
is contractible to a point already in $X$, therefore 
\[
\rank H_3( \widetilde X,\Z) \le
\rank H_3(X_\epsilon,\Z)-2 
\]
\end{remark}

\noindent
In what follows, we will assume that $X\subset \PP^{4}$ is a nodal hypersurface of degree $d$ with $m$ nodes, $l\le m$ of which are 
in general position. Further, let $Y\subset\PP^{4}$ be  a smooth hypersurface of the same degree and 
let $a=\dim H^{0,3}(Y)=\dim H^{3,0}(Y)$ and $b=\dim H^{1,2}(Y)=\dim H^{2,1}(Y)$. Then more is true in (b) of theorem \ref{taut}:

\begin{cor}
 \label{precise}
 If $X$, $\widetilde{X}$ and $Y$ are as before, then 
  \[
 a+b-l = \frac{1}{2} \dim H^3(Y,\C) - l =\dim \Hs_2(Y)-l =  \dim \Hs_2^1\le \dim F^2H^3(\widetilde{X},\C),
 \]
 and therefore  $\rank \, H_3(\widetilde{X},\Z)\ge 2a + 2b - 2l$. 
\\ \\
Since for every point in  a hg set $T$ the rank of $H_3(\widetilde{X},\Z)$ drops off by two with respect to 
the rank of $H_3(Y,\Z)$, the inequality above imposes an upper bound for the number of nodes in any hg set. In particular lemma \ref{taut} shows
that the number of vanishing cycles is at most $l$, i.e.,
there can not be more nodes forming an hg set on $X$ than the number of nodes in general position.
\\ \\
If the $l$ nodes in general position form an hg set on $X$ we actually have
  \[
 \dim \Hs_2^1= \dim \Hs_2(Y)-l = \dim F^2H^3(Y,\C) -l = \frac{1}{2} \dim H^3(Y,\C) - l = a + b -l,
 \]
so that, in this case, theorem  \ref{taut} implies that the map
 $\Hs_2^1\lra F^2H^3(\widetilde{X},\C)$ is in fact an isomorphism.
\end{cor}
The  discussion in corollary \ref{precise}  also implies the following:
\begin{cor}
\label{cor:rank}
If $X\subset\PP^4$ is a nodal hypersurface of degree $d$, where $\Sigma$ 
consist of $m$ nodes in general position and is an hg set on $X$, then
$\dim H^{2,1}(\widetilde{X})=b-m$, $\dim H^3(\widetilde{X},\C) = \rank H_3(\widetilde{X},\Z) = 2a + 2b - 2m$. 
In particular $m\le h^{2,1}(Y)$, 
where $Y\subset \PP^4$ is a smooth hypersurface of degree $d$.
\end{cor}

\begin{remark}
\label{bound}
In particular for a quintic hypersurface $X$ on $\PP^4$ we obtain the nice 
bound $m\le 101$ for the number of nodes in general position which constitute also an hg set. In this case 
(see \cite{BN}, \cite{Ni} and \cite{Vs}) this bound is almost sharp. Observe that the maximal number of
nodes for a quintic hypersurface is expected to lie between 130 and 135, but they do not lie 
in general position (i.e., they impose less than 130 conditions, illustrating one form of the  Cayley-Bacarach theorem). 
If $\Sigma$ is a finite set of nodes, how many independent
conditions does $\Sigma$ impose on homogeneous polynomials of degree $d$ passing through $\Sigma$ and
how many of them form an hg set?
This is not 
the original  formulation of the Cayley-Bacharach theorem but a form of this type of theorem ( see also \cite{egh} p. 297 ).
The exact relation will be given 
by the defect of $X$ considered  in  proposition \ref{defail} of \S \ref{mixed}.
\end{remark}
\noindent
With the same techniques, one can prove a similar result for surfaces on $\PP^3$ and curves on $\PP^2$.
\\
%In the rest of the paper we will assume that nodes in general position are also in hgp.
\\
\\
\section{Generalized Hodge numbers}
\label{general}

Following Danilov and Khovanski\v{\i\ } (see \cite{danilov} \S 1, in particular definition 1.5 proposition 1.8, corollary 1.9 and 1.10), we define the generalized Hodge numbers: 
\begin{eqnarray*}
e^{p,q} = e^{p,q}(X) \stackrel{def} = \sum_k (-1)^k h^{p,q}(H^k_c(X))
\end{eqnarray*}
as well as the generalized Euler characteristic polynomial
\begin{eqnarray*}
e(X;x,\bar{x}) \stackrel{def}= \sum_{p,q} e^{p,q}(X)x^p\bar{x}^q
\end{eqnarray*}
 \\
 which in the sequel we will simply write $e(X)$ and  
 $ \mbox{coeff}_{e(X)}() $ is the coefficient of the term in parenthesis.
 We summarize some well known results about this polynomial (see \cite{danilov}) in a single lemma.
\begin{lem} 
\label{euler} \ \\
\begin{itemize}
\item Suppose $X$ is a disjoint union of a finite number of locally closed subvarieties $X_i, \;\; i\in I$. Then $e(X)=\sum_i e(X_i)$. \\
\item If $f:X\lra Y$ is a bundle with fiber $F$ which is locally trivial in the Zariski topology, then $e(X) = e(Y)\times e(F)$. \\
\item If $X$ is a point, then $e(X)=1$. \\
\item $e(\PP^1)=1+x\bar{x}$. \\
\item $e(\PP^n) = 1 + x\bar{x} + \ldots + (x\bar{x})^n$. \\
\item Let $\pi:\widehat{X}\lra X$ be the blow up of $X$ along a subvariety $Y$ of codimension $r +1$ in $X$. Then
\begin{eqnarray*}
e(\widehat{X})=e(X) + e(Y)[x\bar{x}+ \ldots + (x\bar{x})^r].
\end{eqnarray*}
\end{itemize}
\end{lem}
\ \\ 
As an application of the above lemma we will  compute the generalized Euler polynomial of $X$ for a 
projective hypersurface on $\PP^4$ of degree $d$ with precisely $m$ nodes ($l$ of which are in general position) as the singular locus $\Sigma$. To fix notation, let $\widehat{\PP}^4$ be the blow up of $\PP^4$ along $\Sigma$, $\widehat{X}$ be the inverse image of $X$ on $\widehat{\PP}^4$ and $\widetilde{X}$ be the strict transform of $X$ and $Y$ a non-singular hypersurface of degree $d$ on $\PP^4$. Further, let $\widehat{\Sigma}$ be the inverse image of $\Sigma$ and $\widetilde{\Sigma}=\widehat{\Sigma}\cap \widetilde{X}$. \\ \\
Outside the singular locus the blow up is an isomorphism, therefore one has the following quasi-projective varieties:
\begin{eqnarray*}
X- \Sigma \stackrel{def} = W \cong \widehat{W}\stackrel{def} = \widehat{X} - \widehat{\Sigma}  \cong \widetilde{X} - \widetilde{\Sigma} \stackrel{def} = \widetilde{W} .
\end{eqnarray*}
\\ 
Now, we recall Bott's theorem on the particular situation of $\PP^n$ (\cite{bott} theorems IV and IV'):
 
 \
 \begin{eqnarray} 
 \label{cps} 
 H^p( \PP^n, \Omega^q)= \left\{ \begin{array}{cc}
                                     0    & \mbox{for   $ p \neq q $}, \\
                                     \C    & \mbox{for $ p=q \le n $}
                                 \end{array}
                                     \right.
  \end{eqnarray}
 
  \noindent
 and in particular for $n= 4$:  $ e(\PP^{4})=  1 + x\overline{x} + x^{2}\overline{x}^{2} + x^{3}\overline{x}^{3} + x^{4}\overline{x}^{4}$. 
 It follows inmediately that  $ \Gr_F^j H^n( \PP^4)  = H^{j,n-j}(\PP^4)$  and 
 the only non-zero graded part is when $\mbox{coeff}_{e(\PP^4)}(x^{j}\overline{x}^{n-j}) = 1 $ hence
 \begin{eqnarray}
  \label{grp}
  \Gr_F^2 H^4(\PP^4) = H^{2,2} = \C. 
  \end{eqnarray}
  
 \noindent Also 
 $$
  e(\widehat{\PP}^4) = e( \PP^{4}) + e(\Sigma) ( x \overline{x} +  \ldots + (x \overline{x})^3)
  $$
 using that $ e(\Sigma) = m$ and substituting in the above formula:
 \begin{eqnarray}
 \label{gproj}
 h^{p,q}( \widehat{\PP}^{4}) = \left\{ \begin{array}{cr}
                                        0   &\mbox{if $ p \neq q$},  \\
                                        1   & \mbox{ if $ p = q = 0$},  \\
                                        m+1 & \mbox{if  $ 1\le p = q  \leq 3 $}, \\
                                        1    & \mbox{ if $ p= q= 4$}.
                                        \end{array}
                                        \right.
\end{eqnarray}
It follows that $h^{1,1}(\widehat{\PP}^{4}) = h^{2,2}( \widehat{\PP}^{4}) =  h^{3,3}(\widehat{\PP}^{4}) = m+1 $.
After these basic preliminaries, the main result in this section is:
\begin{prop} 
\label{ghn} Let $X, \widetilde{X}$ and $Y$ as above, then 
\[
e(\widetilde{X}) = 1 + (m+1) x \overline{x} - ax^{3} - (b-l)x^{2}\overline{x}  - (b-l)x \overline{x}^{2} - a\overline{x}^{3} + (1+m) x^{2} \overline{x}^{2} + x^{3}\overline{x}^{3}
\]
and
\[
e( X) = 1+(1-m) x \overline{x}- ax^{3}- (b-l)x^{2}\overline{x} - (b-l)x \overline{x}^{2} - a\overline{x}^{3} + x^{2}\overline{x}^{2}  + x^{3} \overline{x}^{3}, \hskip2cm{} 
\]
where $a=h^{3,0}(Y)$, $b=h^{2,1}(Y)$.
\end{prop}
\pf 
Observe that $ \widehat{\Sigma} = \cup_{ x \in \Sigma} E_{x} $ and by  cutting each $E_{x} $ with $ \widetilde{X}$ we obtain a quadric surface $Q_{x}$ hence 
$ e( \widetilde{\Sigma}) = \Sigma_{x} e(Q_{x}) $ but each summand is equal to 
\[ e(\PP^1\times\PP^1)=e(\PP^{1})^{2}= 1 + 2x\overline{x} + x^{2} \overline{x}^{2},\]
so
\begin{eqnarray}
\label{cua}
 e( \widetilde{\Sigma}) = m\left( 1 + 2x\overline{x} + x^{2} \overline{x}^{2} \right).
\end{eqnarray}
 Moreover, $ e^{p,q}(\widetilde{W}) = e^{p,q}(\widetilde{X}) -e^{p,q}(\widetilde{\Sigma})$ and 
 \[ e^{p,q}(X) = e^{p,q}(W) + e^{p,q}(\Sigma) = e^{p,q}(\widetilde{W}) + e^{p,q}(\Sigma) = e^{p,q}(\widetilde{X}) - e^{p,q}(\widetilde{\Sigma}) + e^{p,q}(\Sigma).\] Since $h^3(\widetilde{X})= 2a + 2b - 2l $ (see corollary \ref{cor:rank}), Lefschetz Hyperplane Theorem tell us that:
\begin{eqnarray}
\label{ghn2}
\end{eqnarray}
 \begin{eqnarray*} 
  e(\widetilde{X}) = 1 + (m+1) x \overline{x} - ax^{3} - (b-l)x^{2}\overline{x}  - (b-l)x \overline{x}^{2} - a\overline{x^{3}} + (1+m) x^{2} \overline{x}^{2} + x^{3}\overline{x}^{3}. 
  \end{eqnarray*}
  Finally, $ e(X) = e(\widetilde{X}) -(m + 2m x \overline{x} + mx^{2}\overline{x}^{2}) + m$. The result follows directly by substituting the value of  $e(\widetilde{X})$ in  equation (\ref{ghn2}).
\hfill \qed
\ \\  \ \\
Using the  Hodge numbers of the total transform $\widehat{\PP}^4$  given by equation (\ref{gproj}) we can conclude:

\begin{cor}
\label{def}
 In lemma \ref{taut}  of  \S \ref{hyp} in diagram (\ref{ctaut}): $ d( A^3_1(X)) = \widetilde{d}( A^3_1(\tilde X)) = 0 $  hence $\Cok(d) = \Hs^1_2 =A^4_2(X, \Sigma )$  and 
 $\Cok(\widetilde{d}) = \Hs_2(\widetilde{X})  = A^4_2(\widetilde{X})$.
\end{cor}
\pf  Since  $ \tilde X $ is smooth, then the hodge numbers $e^{p,q}(\widetilde{X})=(-1)^{p+q}h^{p,q}(\widetilde{X})$, in particular $ h^{2,0} = h^{0,2} = 0 $  by the computation above. This implies that
$ H^0( \tilde X,  \widehat{\Omega} _{\tilde X}^2 ) \subset 
H^0( \tilde{X},\Omega _{\tilde X}^2 ) = 0 $. Recall the exact sequence of residues in  \cite{griffiths2} Lemma 10.9  ii): 
\[
 0 \rightarrow  \widehat{\Omega}_{\tilde{\PP}}^q  \rightarrow  \widehat{\Omega}_{ \tilde{\PP}}^q(1)  \rightarrow 
 \widehat{\Omega}_{ \tilde X}^{q-1} \rightarrow 0
\]
and its associated long sequence for $ q= 3$:
\[
0 \rightarrow H^0(  \tilde{\PP}, \widehat{\Omega}_{ \tilde{\PP}}^3 ) \rightarrow  H^0(\tilde{\PP} ,\widehat{\Omega}_{ \tilde{\PP}}^3(1) ) \rightarrow 
H^0( \tilde X,  \widehat{\Omega} _{\tilde X}^2 )  \rightarrow \cdots
\]
also  $H^0(  \tilde{\PP}, \widehat{\Omega}_{ \tilde{\PP}}^3 ) \subset H^0(  \tilde{\PP}, \Omega_{ \tilde{\PP}}^3 ) = 0 $ (see equation (\ref{gproj}) above ) 
since the last term for the  above sequence 
 is  already zero  so must be  the middle term. In particular, $ d( A^3_1(X)) \subset \widetilde{d}( A^3_1(\tilde X)) = 0 $.
\hfill \qed

\section{Mixed Hodge structure of a nodal 3-fold}
\label{mixed}

Given a singular scheme $X$ defined over $\C$, Guillen, Navarro \textit{et. al.} defined a \emph{cubical hyperresolution} $X_{\bullet}$ of $X$ (see\cite{Navarro}, Expos\'e III, proposition 3.3) which induces a spectral sequence
\[
E_1^{p,q}=H^q(X_p,\C)\Rightarrow H^{p+q}(X,\C)
\]
providing a natural Mixed Hodge Structure on $H^{p+q}(X,\C)$ (we set
$X_p \stackrel{def} =  \sqcup_{|\alpha|=p+1} X_{\alpha}$).
\\ \\
In our situation, a cubical hyperresolution can be constructed from the following pullback diagram
\[
\xymatrix{
\Sigma\times_X\widetilde{X}\ar[d]\ar@<-0.25ex>[r] & \widetilde{X}\ar[d]^{\pi}\cr
\Sigma\ar[r] & X
}
\]
Since $\Sigma\times_X\widetilde{X}\cong\widetilde{\Sigma}$, the projection to the first factor gets identified with $\pi_|$, the restriction of $\pi$ to $\widetilde{\Sigma}$, while the projection to the second factor gets identified with the natural inclusion $i:\widetilde{\Sigma}\hookrightarrow\widetilde{X}$, yielding the cubical hyperresolution
\[
\xymatrix{
X_1\ar@<0.1ex>[r]^{i}\ar@<-0.5ex>[r]_{\pi_|}& X_0\ar[r] & X,
}
\]
where $X_1=\widetilde{\Sigma}$ and $X_0=  \widetilde{X}\sqcup\Sigma$. Therefore $E_1^{0,q}=H^q(X_0,\C)$, $E_1^{1,q}=H^q(X_1,\C)$ and $E_1^{p,q}=0$ for all $p\ge 2$. Clearly this spectral sequence degenerates at $E_2$, so we have
\[
0\to E_2^{1,2}\to H^3(X,\C) \to E_2^{0,3} \to 0,
\]
where \[E_2^{0,3}=\Ker (\xymatrix{H^3(X_0,\C) \ar[r]^{\pi_|^*-i^*} & H^3(X_1, \C)})\] and 
\[E_2^{1,2} = H^2(X_1,\C)/(\im(\xymatrix{H^2(X_0,\C)\ar[r]^{\pi_|^*-i^*} & H^2(X_1,\C)}).\]
Since $H^3(X_0,\C) = H^3(\widetilde{X},\C)$ and  $H^3(X_1, \C)=0$, then $E_2^{0,3}= H^3(\widetilde X,\C)$ is a pure Hodge structure of weight 3. \\ \\ 
Similarly, 
$H^2(X_0,\C) = H^2(\widetilde{X},\C)\cong\C^m$ and $H^2(X_1,\C)= H^2(\widetilde{\Sigma},\C) \cong \C^{2m}$, so
$E_2^{1,2}\cong\C^m$ is a pure Hodge structure of weight 2 and we recover the Clemens-Schmidt exact sequence
\begin{eqnarray} 
\label{Clemens-Schmidt}
0\lra W_2H^3(X,\C)\to H^3(X,\C)\to H^3(\widetilde{X},C)\to 0
\end{eqnarray}
with $W_2H^3(X,\C)= E_2^{1,2}\cong\C^m$, which is to be expected for a cubical hyperresolution, as pointed out in \cite{peters}, Corollary 5.42.
\\ \\
Remember that, in virtue of  theorem \ref{ghn}, if $\Sigma$ consists of $m$ nodes, where precisely $l$ of them are in general position (and assuming they are also in homologically good position), one has
\begin{enumerate}
\item $\dim H^3(\widetilde{X}, \C) = 2a+2b-2l$ \\
\item $H^3(\widetilde{X}, \C) \cong \oplus \widetilde{H}^{i,j}$, 
\ \\ \\
where $\dim \widetilde{H}^{0,3} = \dim \widetilde{H}^{3,0} = a$ and $\dim \widetilde{H}^{1,2}=\dim \widetilde{H}^{2,1} = b-l$. \\
\item $\dim H^3(X,\C) = 2a+2b-2l + m$.
\end{enumerate}
\noindent
Moreover, in this situation we have:

\begin{prop}
\label{GrF}
\[
{\Gr_F}^{k}H^3(X,\C)  = \left\{ \begin{array}{cr}
                                       \C^a   \hfill      &\mbox{if $ k = 0,3$},  \\
                                        \C^{b-l+m}   & \mbox{ if $ k = 1$}, \\
                                       \C^{b-l}  \hfill & \mbox{ if $ k = 2$ }.
                                         \end{array}
                                        \right.
\]
Observe that $m-l$ is precisely the failure of $\Sigma$ to impose independent conditions on homogeneous polynomials of degree 5 (see remark \ref{bound}). 
\end{prop}

\begin{ej}
\label{ej:quintic}
Let $X$ be the quintic threefold on $\PP^5$ defined by the equations $p_1 = 0$ and $4 p_5 -15 p_2 p_3 = 0$, where $p_k=\sum_0^5 x_i^k$ is the $k$-th power symmetric function. Then the singular locus of $X$ consists of precisely $100$ nodes which are the
orbits of $(1:-1:1:-1:1:-1)$ and of $(1:-1:1:-1:z:-z)$ under the symmetric group on six letters $S_6$, where 
$7z^2+16=0$. For this quintic threefold,  using the  kernel extension { \tt PLURAL} of {\tt SINGULAR} 2-0-6 (see \cite{plural} ) we have written a program that
allows us to conclude that the 100 nodes impose only 86 conditions on the space of quintics passing through them, so in this case $l=86<100=m$.
It is not difficult to see that this quintic threefold is actually a singular Calabi-Yau threefold on $H=V(p_1)\cong\PP^4$. 
\\ 
As Candelas, de la Ossa \'et al. have shown in \cite{candelas} , $\dim H^3(Y,\C) = 204$ for a smooth quintic threefold on $\PP^4$,
and if additionally the nodes in general position form an hg set, then  corollary (\ref{GrF}) can be written as:
\[
{\Gr_F}^{k}H^3(X,\C)  = \left\{ \begin{array}{cr}
                                       \C   \hfill      &\mbox{if $ k = 0,3$},  \\
                                        \C^{115}   & \mbox{ if $ k = 1$}, \\
                                        \C^{15}  \hfill & \mbox{ if $ k = 2$}. 
                                        \end{array}
                                        \right.
\]
and $\dim W_2 H^3(X,\C) = 100$.
\end{ej}

Recall  that we have a  commutative diagram of long exact sequences with compact support:
\[
\xymatrix{
\dots \ar[r] & H^i_c(U)\ar[r]\ar[d]^{\pi^*}_{\cong} & H^i(X)\ar[r]\ar[d]^{\pi^*} & H^i(\Sigma) \ar[r]\ar[d]^{\pi^*} &
H^{i+1}_c(U)\ar[r]\ar[d]^{\pi^*}_{\cong} & \dots\cr
\dots \ar[r] & H^i_c(\widetilde{U})\ar[r] & H^i(\widetilde{X})\ar[r] & H^i(\widetilde{\Sigma}) \ar[r] &
H^{i+1}_c(\widetilde{U})\ar[r] & \dots
}
\]
Since $\Sigma$ is zero dimensional, then $H^i(\Sigma)=0$ for all $i>0$. In particular:
\[
\xymatrix{
0 \ar[r] & H^4_c(U)\ar[r]^{\cong}\ar[d]^{\pi^*}_{\cong} & H^4(X)\ar[r]\ar[d]^{\pi^*} & 0 \ar[r]\ar[d]^{\pi^*} &
H^{5}_c(U)\ar[r]^{\cong}\ar[d]^{\pi^*}_{\cong} & H^5(X)\ar[r]\ar[d]^{\pi^*} & 0\cr
0 \ar[r] & H^4_c(\widetilde{U})\ar[r] & H^4(\widetilde{X})\ar[r] & H^4(\widetilde{\Sigma}) \ar[r] &
H^{5}_c(\widetilde{U})\ar[r] &H^5(\widetilde{X})\ar[r] & 0
}
\]
is exact and commutative. From the generalized Hodge numbers, equations (\ref{cua}) and (\ref{ghn2}), we have: $H^i(\widetilde{\Sigma})=0$ for $i>4$, 
 $H^3(\widetilde{\Sigma})=0$, $H^4(\widetilde{\Sigma})\cong\C^m$, 
 $H^4(\widetilde{X})\cong\C^{m+1}$   and  $H^5(\widetilde{X}) = 0$. Therefore the second row of the above diagram simplifies to:
\[
 0 \rightarrow \C^{\beta}  \rightarrow \C^{m+1} \rightarrow \C^m \rightarrow  \C^r \rightarrow  0  
\]
where $\beta$ is the fourth Betti number of $X$.
Applying  the Euler characteristic to this exact sequence: $ \beta-(m+1) + m -r =0 $ hence $ r = \beta-1$.
It follows that $H^5(X)\cong H^5_c(\widetilde{U}) \cong  \C^{\beta-1}$ for some $\beta \le m+1$.
\\ \\
On the other hand, Clemens in  \cite{Cle} and  later Werner in \cite{We} have introduced the following  Mayer-Vietoris type  exact sequence:
\[
\xymatrix{
0\ar[r] & H_4(Y)\ar[r] & H_4(X)\ar[r]^k&\mathcal{R}\ar[r]^{b} & H_3(Y)\ar[r]^{\gamma} & H_3(X)\ar[r] & 0
}
\]
 where $Y$ is a smooth  threefold  of the same degree as $X$ and $\mathcal{R}$ is a free $\Z$-module of rank $m=|\Sigma|$. This  allows us to compute the 
 defect of $X$ as $ \delta \stackrel{def}{=} \rank(  \im(k) ) $. As a consequence of their definition 
 they show that $ \beta_2(X) = 1$  and $\delta = \beta -1 $.
 \begin{cor}
  If the  $l$  double points are in  general position (resp. form a hg set)  then $ \delta  \geq 2(m-l)$ 
  (resp.  $ \delta  =  2(m-l)$). In particular, 
  $ l \geq  \frac{m}{2} $. 
 \end{cor}
 \pf 
 If the $l$ points are in  (algebraic) general position  then $ \dim H^3( \tilde X)  \geq  2(a +b -l) $ and by  equation 
 (\ref{Clemens-Schmidt}): $ \dim H^3(X)  =  m + \dim H^3( \tilde X) \geq 2(a +b -l) + m $, but 
 $ \rank(\im(b))  = h^3(Y) - h^3(X)  = 2(a +b) -h^3(X) \leq  2(a+b) -2(a+b -l) -m = 2l -m $ 
 ( inequality is an equality if all the double points form a hg set). Hence $ \delta=\rank(\Ker(k))= m  - \rank(\im(b))\ge m - (2l-m) = 2(m-l)$  and 
 $m+1\ge \beta=\delta+1\ge 2(m-l)+1$. Therefore $m\le 2l$. \hfill \qed
 \\ \\
 In order to find an exact relation  between $ \delta$ and the failure of $\Sigma$  to impose   linearly independent conditions  on polinomials of degre $2d-5$ 
  (compare with  remark \ref{bound} of \S \ref{hyp})  we shall use and prove the following:
 \begin{prop}
 \label{defail}
  If the l double points are in (algebraic) general position and form an hg set, then
  \[
\delta = m- l + a + b - \binom{ 2d-1}{4}.
\]
 \end{prop}
 \pf
By \cite{We} (see Satz Kap. IV p.27) $ \delta = m-\binom{2d-1}{4} + \dim( A^4_2(X, \Sigma)) $. By corollary \ref{def} the last term
$ A^4_2(X, \Sigma)  = \Hs^1_2$ and  by the assumption on $\Sigma $ the dimension of the latter is equal to $a+b-l$. \hfill
\qed 

\begin{remark} 
\label{defail2}
The significance of the last corollary is that the difference between the 
defect and the failure of $\Sigma$ to impose conditions on polynomials of degree  $2d-5$ is  equal to  $ a + b  - \binom{ 2d-1}{4} $ which depends only on the degree of $X$ and the dimensions
$ h^{3,0}, h^{2,1} $  of a smooth $Y$ of the same degree as $X$.  
\end{remark}

\section{Equisingular families}
\label{equi}

Let 
\begin{eqnarray*}
\xymatrix{
\bar{\sX} \ar@{^{(}->}[r]\ar[rd]_{f} & \PP^4 \times\PP^1\ar[d]\cr                                  
                                                 & \PP^1
}
\end{eqnarray*}
be a Lefschetz pencil of hypersurfaces on $\PP^4$, where the vertical arrow is the projection on the second factor, and assume that there is a maximal
non empty open subset $B\subset\PP^1$ over which the family 
\begin{eqnarray*}
\xymatrix{
\sX = f^{-1}(B) \ar@{^{(}->}[r] \ar[d]_{f}& \bar{\sX}\ar[d]^{\tilde{f}}\cr
 B  \ar@{^{(}->}[r]  & \PP^1
}
\end{eqnarray*}
is real analitically trivial and such that the singular locus  $\Sigma_t$ of every fiber $X_t$ consists 
of exactly $m$ nodes. Then the higher direct image $\H=R^3f_*\C$ is a local system, with fiber $H^3(X_t,\C)$ admiting a MHS. 
%%%%%%%%%%%%%%%%%%%%%%%%%%%%%%%%%%%%%%%%%%%%%%%%%%%%%%%%%%%%%%%%%%%%%%%%%%%%%%%%%%%%%%%%%%%%%%%%%%%%%%%%%%%%%%%%%%%%%%%%%%%%%%%%%%%%%%%%%
%which then becomes a VMHS.
\\ \\
For a fixed $t\in B$, let $\widehat{\PP}^4$ be the blow up of $\PP^4$ along $\Sigma_t$ and $\widetilde{X}_t$ be the strict transform of $X_t$. Further, let $\widehat{\Sigma}_t$ be the inverse image of $\Sigma_t$ (i.e., the disjoint union of the exceptional divisors along the $m$ nodes) and $\widetilde{\Sigma}_t = \widehat{\Sigma}_t\cap \widetilde{X}_t$. Since the multiplicity of every point in $\Sigma_t$ is 2, then $\widetilde{X}_t$ is a projective, non singular variety and we have a diagram
\begin{eqnarray*}
\xymatrix{ 
            {\Sigma_t} \ar@{^{(}->}[r]   &   {X_t} \ar@{^{(}->}[r]   &  {\PP}_4  \\
           \widetilde{\Sigma}_t \ar@{->}[u]            \ar@{^{(}->}[r] &  \ar@{->}[r]  \widetilde{X}_t \ar@{->}[u]_{{\pi}}  \ar@{^{(}->}[r]  &  \widehat{{\PP}}_4  \ar@{->}[u]^{ \pi }  }       
\end{eqnarray*}

Let $\widetilde{\sX}\stackrel {\tilde{f}} \longrightarrow  B$ be the smooth family formed by the union of 
$\widetilde{X}_t$ along $B$ (also a Lefschetz pencil). Then the 
higher direct image $\widetilde{\H}=R^3\widetilde{f}_*\C$ is a VHS on $B$, in particular we have a  GM
-connection \[\widetilde{\nabla}^{GM}:\widetilde{\sH}^3\to \widetilde{\sH}^3\otimes\Dif_B^1,\]
where $\widetilde{\sH}^3\stackrel{def}{=}R^3\tilde{f}_*\C\otimes\sO_B$. As seen in the introduction, at every point $p\in\PP^1\backslash B$ there exist a limit MHS $(W_{\bullet}, F{\bullet})$, as well as an extension of the monodromy operator $T$, inducing the weight filtration $W_{\bullet}$, and such that the corresponding nilpotent operator $N=\text{log}(T_u)$ has nilpotence degree $\le 4$. 
\\ \\
 \begin{ej} 
 \label{non-singular}
It is not difficult to see, using the above description and the notation 
on proposition (\ref{GrF}), that for a smooth quintic threefold $X\subset\PP^4$ one has $a=1$ and $b=101$ (see \cite{candelas}).
Moreover, if $X_t$ is the smooth family of quintic threefolds in $\PP^4$ given by $x^5+y^5+z^5+w^5+u^5-5txyzwu$, it has been shown by Candelas \'et . al.  that the GM-connection induces
a maximal unipotent map  on $H^3(X_t,\C)$, for any $t$, whose nilpotent part $N$ satisfies
$N(H^{p,3-p})\subset H^{p-1,3-p+1}$ for $0\le p\le 3$ with $N^3\ne 0$ but $N^4=0$. In particular 
one has an splitting of the Hodge structure:
\begin{eqnarray*}
H^3(X_t,\C) = J\oplus_{i=1}^{100} V_i(-1)
\end{eqnarray*}
where $J$ is a Hodge structure of weight 3 and type $(1,1,1,1)$ and each
$V_i(-1)$ is a Hodge structure of weight 3 and type $(0,1,1,0)$, associated to a Hodge structure $V_i$ of  weight 
one and type $(1,1)$ (see  also \cite{candelas} for the quintic family  of threefolds
in connection with mirror symmetry). Here, as usual, $V_i(-1)=V_i\otimes \Z(-1)$ and $\Z(-1)$ is the Tate-Hodge structure of weight 2.
\end{ej}

\begin{ej}
  \label{jord}
More generally, for a pencil of Calabi-Yau threefolds on $\PP^4$ we have $$\dim (\widetilde{H}^{0,3}) =  \dim (\widetilde{H}^{3,0})= 1 $$ and  $ k = \dim (\widetilde{H}^{1,2}) =
   \dim(\widetilde{H}^{2,1}) $  hence $\widetilde{H} \cong \C^{2k+2}$.
   \\ \\
   In the same spirit as the example given in \cite{candelas} and example \ref{ej:quintic}, keeping the notation there, 
  for the case $ n= 5$, consider the pencil of quintic hypersurfaces in $ H= \PP^4$ defined by \\ \\
\centerline{
 $
 f_{(\alpha, \beta)} = \alpha p_5 - \frac{5(\alpha + \beta)}{6} p_2 p_3.
$}
\\ \\
Let $ {\sM} \subset H  \times \PP^1$ be the corresponding incidence family. Clearly, for each $ (\alpha :\beta)$ we 
have a  quintic $ {\sM}_{(\alpha:\beta)} \subset  \PP^{4}$. 
  This family has  already been introduced and studied  by Van Straten in \cite{Vs}. In {\it loc.cit }
  ( see Theorem 2 ), he shows  that ${\sM}_{(\alpha:\beta)}$ is a singular variety for a general value of  $( \alpha: \beta) =  ( \frac{\alpha}{\beta}:1 ) $,~except for the quintics associated to:
 \[
  q_1=25, q_2=1, q_3=-3, q_4=0, q_5=-2, q_6 = \infty .
  \] 
  For $t\in\PP^1-\{q_1,\dots , q_6\}$, the singular  locus, $\Sigma_{t}=\Sing(\sM_{t})$ consist of  
  $100$ nodes (compare with the bound $m\le 101$ computed in remark \ref{bound}). 
  In example \ref{ej:quintic}, we have seen that only $86$ of these nodes are in general position, therefore 
  for a general member of this family we have $ \dim H^3(X_t)=132 $ and
\[
Gr^k_FH^3(X_t,\C)  \cong \left\{ \begin{array}{cr}
                                        \C   \hfill      &\mbox{if $ k = 0,3$},  \\
                                        \C^{115}   & \mbox{ if $ k = 1$}, \\
                                        \C^{15}  \hfill & \mbox{ if $ k = 2$}.
                                        \end{array}
                                        \right.
\]
\noindent
while $\dim H^3(\widetilde{X}_t)=32 $ and
\[
\rank \widetilde{H}^{k,3-k} = \left\{ \begin{array}{cr}
                                       1   \hfill      & \mbox{if $ k = 0,3$},  \\
                                       15   & \mbox{ if $ k = 1,2 $}.
                                        \end{array}
                                        \right.
\]
 \end{ej}
%%%%%%%%%%%%%%%%%%%%%%%%%%%%%%%%%%%%%%%%%%%%%%%%%%%%%%lmhs%%%%%%%%%%%%%%%%%%%%%%%%%%%%%%%%%%%%%%%%%%%%%%%%%%%%%%%%%%%%%%%%%%
Before we study the LMHS  of the VPHS given by  $\widetilde \sX \rightarrow B$ we  introduce a very well known 
inductive method to calculate the monodromy weight filtration and advise  the reader interested in the main result
 to skip to proposition \ref{lmhs}. For that  let $m$ be an integer, $ H_{\Q}$ be a $\Q$-vector 
 space and $N:H_{\Q}\to H_{\Q}$ be a nilpotent endomorphism such that $ N^{m+1} = 0 $ but $N^m\ne 0$.
Following Donagi (see \cite{griffiths3}, remark on page 69), we can introduce the following  
$\Q$-spaces for $r,s $ positive integers satisfying $ r \le m, s \leq m+1 $:  
$ M_{r,s} = \Ker N^{m-r} \cap \im N^s  $. These spaces satisfy the following relations: \\
$ M_{0, s}  \supset  M_{1,s} \supset \ldots \supset M_{m,s} = 0 $ and 
similarly  $M_{r,0} \supset M_{r,1}  \supset  \ldots \supset M_{r,m+1} =0 $. Observe that the nilpotent operator $N$ admits a natural extension to $N:H_{\C}\to H_{\C}$.
\\ \\
Consider an increasing filtration on $H_{\Q}$ given by the $\Q$-vector spaces  : \[  W_q \stackrel{def}{=} < \Sigma_{ 2m-q- 1 =r +s}  M_{r,s} >,\] for $0  \leq q  \leq 2m-1 $, while $W_{2m}\stackrel{def}{=} H_{\Q}$.
 
\begin{ej}
\label{app}
Since the sum in the formula is internal  one need not compute all terms in the formula above and many redundant terms occur. If one represents the lattice of subspaces $M_{r,s}$ as integral points in the plane, observe first that $M_{m,j}=0$ and $M_{j,m+1}=0$ for all $j$, since $N^0=id$ (the identity) and $N^{m+1}=0$; therefore the relevant terms lie in the integral points of a finite array of  $(m-1)\times m$. Moreover, if $s-r\ge 1$, then $N^{m-r+s}=0$, i.e., $\im N^s\subset\Ker N^{m-r}$ and $M_{r,s}=\im N^s$.
  \\ \\
  The subspace $W_q$ is the sum of the subspaces represented by the integral points lying on the line $ L_q  \stackrel{def}{=} \{ (r, s) | r+s = 2m-q-1 \} $ and we want to know which of the corresponding subspaces $M_{r,s}$ actually contribute to the sum. 
  \\ \\
  \[
  \xy <2cm,0cm>:
  (0,0) *=0{\bullet}="{\bullet}" ;
  (0,1) *=0{\bullet}="{\bullet}" ;
  (0,2) *=0{\bullet}="{\bullet}" ;
  (0,3) *=0{\bullet}="{\bullet}" ;
  (0,4) *=0{\bullet}="{\bullet}" ;
  (0,5) *=0{\bullet}="{\bullet}" ;
  (1,0) *=0{\bullet}="{\bullet}" ;
  (1,1) *=0{\bullet}="{\bullet}" ;
  (1,2) *=0{\bullet}="{\bullet}" ;
  (1,3) *=0{\bullet}="{\bullet}" ;
  (1,4) *=0{\bullet}="{\bullet}" ;
  (1,5) *=0{\bullet}="{\bullet}" ;
  (2,0) *=0{\bullet}="{\bullet}" ;
  (2,1) *=0{\bullet}="{\bullet}" ;
  (2,2) *=0{\bullet}="{\bullet}" ;
  (2,3) *=0{\bullet}="{\bullet}" ;
  (2,4) *=0{\bullet}="{\bullet}" ;
  (2,5) *=0{\bullet}="{\bullet}" ;
  (3,0) *=0{\bullet}="{\bullet}" ;
  (3,1) *=0{\bullet}="{\bullet}" ;
  (3,2) *=0{\bullet}="{\bullet}" ;
  (3,3) *=0{\bullet}="{\bullet}" ;
  (3,4) *=0{\bullet}="{\bullet}" ;
  (3,5) *=0{\bullet}="{\bullet}" ;
  (4,0) *=0{\bullet}="{\bullet}" ;
  (4,1) *=0{\bullet}="{\bullet}" ;
  (4,2) *=0{\bullet}="{\bullet}" ;
  (4,3) *=0{\bullet}="{\bullet}" ;
  (4,4) *=0{\bullet}="{\bullet}" ;
  (4,5) *=0{\bullet}="{\bullet}" ;
  (5,0) *=0{\bullet}="{\bullet}" ;
  (5,1) *=0{\bullet}="{\bullet}" ;
  (5,2) *=0{\bullet}="{\bullet}" ;
  (5,3) *=0{\bullet}="{\bullet}" ;
  (5,4) *=0{\bullet}="{\bullet}" ;
  (5,5) *=0{\bullet}="{\bullet}" ;
  (0,0) *+{\bullet}; (5,5)*+{\bullet} **@{-} 
 ? *{\hskip-1cm{\Delta}};
 (0,2) *+{\bullet}; (2,0)*+{\bullet} **@{-} 
 ? *{\hskip1cm{L_7}} ;
 (1,5) *+{\bullet}; (5,1)*+{\bullet} **@{-} 
 ? *{\hskip1cm{L_2}} ;
 (0,0) *+{\bullet}; (0,5)*+{\bullet} **@{-} ;
 (0,0) *+{\bullet}; (5,0)*+{\bullet} **@{-} ;
 (1,4) *+{\bullet}; (2,4)*+{\bullet} **@{ } ;
 ? *{\hskip-0.5cm{\text{}^{Im{N^4}}}};
 (3,4) *+{\bullet}; (4,4)*+{\bullet} **@{ } ;
 ? *{\hskip-0.5cm{\text{}^{Im{N^4}}}} ;
 (0,0) *+{\bullet}; (0,1)*+{1} **@{} 
  *+/v(0,1)/{\hskip-0.5cm{\text{}_{1}}};
(0,1) *+{\bullet}; (0,2)*+{\bullet} **@{} 
  *{\hskip-0.5cm{\text{}_{2}}};
(0,2) *+{\bullet}; (0,3)*+{\bullet} **@{} 
  *{\hskip-0.5cm{\text{}_{3}}};
(0,3) *+{\bullet}; (0,4)*+{\bullet} **@{} 
  *{\hskip-0.5cm{\text{}_{4}}};
 (0,4) *+{\bullet}; (0,5)*+{\bullet} **@{} 
  *{\hskip-0.5cm{\text{}_{5}}};
  (0,5) *+{\bullet}; (0,0)*+{\bullet} **@{} 
 *{\hskip-0.5cm{\text{}_{0}}};
 (0,0) *+{\bullet}; (1,0)*+{\bullet} **@{} 
  *{\hskip-0.5cm{\text{}_{1}}};
(0,1) *+{\bullet}; (2,0)*+{\bullet} **@{} 
  *{\hskip-0.5cm{\text{}_{2}}};
(0,2) *+{\bullet}; (3,0)*+{\bullet} **@{} 
  *{\hskip-0.5cm{\text{}_{3}}};
(0,3) *+{\bullet}; (4,0)*+{\bullet} **@{} 
  *{\hskip-0.5cm{\text{}_{4}}};
 (0,4) *+{\bullet}; (5,0)*+{\bullet} **@{} 
  *{\hskip-0.5cm{\text{}_{5}}};
  \endxy
  \]
 \  \\  \ \\
  If $q=2b$ for some integer $b$, then the intersection of the diagonal $\Delta \stackrel{def}{=} \{ (r, s) | r=s\} $ and $L_q$ is not an integral point, but  $(m-b-1,m-b)\in L_q$. As observed above, since $m-b - (m-b-1)=1$, $M_{m-b-1,m-b}=\im N^{m-b}$. Moreover, all points $(r,s)\in L_q$ lying above the diagonal $\Delta$ satisfy $s\ge m-b>m-b-1\ge r$, therefore $s-r\ge 1$ and
  $M_{r,s}=\im N^s\subset \im N^{m-b}= M_{m-b-1,m-b}$ for all such points $(r,s)\in L_q$, in particular the corresponding $M_{r,s}$ do not contribute anything new to $W_q$.
 \\ \\ 
 If $q$ is odd, then $2m-q-1$ is even and $L_q$ intersects the diagonal $\Delta$ at the integral point $(m - (q+1)/2, m - (q+1)/2)$. In this case, all points $(r,s)\in L_q$ above the diagonal $\Delta$ satisfy $s>m-(q+1)/2>r$ and again $s-r\ge 1$, so $M_{r,s}=\im N^s\subset \im N^{m-(q+1)/2}$. Since $s+(q+1)/2>m$ for all such points, one has 
 $N^{(q+1)/2+s}=0$, i.e., $\im N^s\subset\ker N^{(q+1)/2}$ as well, therefore  $M_{r,s}\subset M_{m-(q+1)/2,m-(q+1)/2}$ and these $M_{r,s}$ do not contribute anything new to $W_q$.
 \\
We can summmarize these calculations saying that the only subspaces $W_{r,s}$ that contribute something new to $W_q$ are those subspaces for which
 $(r,s)\in L_q$ lies below $\Delta$, on $L_q\cap\Delta$ (if $q$ is odd) or inmediatly above $\Delta$ (if $q$ is even). For instance for $m=3$ one has 
  \begin{eqnarray*}
  \begin{matrix}
   W_0 = < N^{3,0} >,\hfill & W_1 = <N^{2,2} >, \hfill & W_2 = < N^{1,2}  + N^{2,1} >, \\
   W_3 = < N^{0.2} + N^{1,1} > , &  W_4 =  < N^{1,0} + N^{0,1} > , & W_5 = < N^{0,0} > .\hfill
   \end{matrix}
  \end{eqnarray*}
As one can see, at most $9$ different summands contribute to  $W_q$'s in this case. More generally, denoting by  $N_m $ the maximal number of \it{different} summands  contributing to the $W_q$'s, an elementary  counting argument shows that for arbitrary $m$ this number equals the number of points $I_m$  in  the isosceles triangle  of height $m-1$ and base $m-1$  bounded below by the diagonal, plus the number of odd integers in the set $\{1,2, \ldots, 2m-3\}$ plus 
  one, giving $ N_m = \frac{m(m+3)}{2} $. In particular for $ m \geq 1$
  $$
  I_m = \frac{(m-1)^2}{2}  \leq  N_m  \leq m^2.
  $$
 This is already true for $ m = 3$ as seen above. Similarly, $N_4 = 14, N_5 = 20 $. 
\end{ej}
\begin{lem}\label{lem:monodromyweight}
The filtration defined above satisfies Morrison's characterization of the weight filtration on $H_{\Q}$ associated to $N$.
\footnote{The formula for $W_k$ given in \cite{griffiths3}, page 69  is incomplete. The procedure  there described is correct, 
however the formula  has a misprint. Here we include a more accurate formula in both cases for lack of another suitable reference.} 
(see \cite{griffiths3}, pages 106 - 107):
\begin{enumerate}
\item $N(W_k)\subset W_{k-2}$,
\item $W_{m-t}/W_{m-t-1}=\im(N^t|_{W_{m+t}/W_{m-t-1}})$,
\item $W_{m+t-1}/W_{m-t-1}=\Ker(N^t|_{W_{m+t}/W_{m-t-1}})$.
\end{enumerate}
\end{lem}

\pf 
It is helpful to visualize the action of $N^t$ on the Lattice formed by the $M_{r,s}$ as follows:
\[
\xymatrix @=0.1cm{
M_{j-1,m-t-j}\ar@{->>}[rrdddd]^{N^t} & \supset \cdots  \supset & M_{t+j-1,m-t-j}\cr
 \cup &  & \cup\cr
\vdots &  & \vdots\cr
\cup &  & \cup\cr
M_{j-1,m-j} & \supset \cdots  \supset & M_{t+j-1, m-j}
}
\]
\\
1) Obviously {\small $N(M_{a,b})=N(\Ker N^{m-a}\cap\im \,N^b)\subset \Ker N^{m-(a+1)}\cap  \im \,N^{b+1}=M_{a+1,b+1}.$}
\\ \\
2) Observe that $W_{m-t} = M_{t-1,m}+\dots + M_{m, t-1}$, since $N^{m+1}\equiv 0$ and $\Ker N^0=0$. 
We claim that $N^t(M_{r,s})=M_{r+t,s+t}$.  Indeed,  the first inclusion is the content of the proof above for 1). For the equality, let  $ x \in \Ker N^{ m-r-t} $ and   $x =\im( N^{s+t}(y))$. Let $z=N^{s}(y)$, therefore $N^t(z) = N^{s+t}(y) = x$ and  $ 0 = N^{m-r-t}(x) =N^{m-r-t}( N^{s+t}(y))= N^{m-r}(N^s(y)) = N^{m-r}(z)$, i.e. $z \in M_{r,s} $. It follows that $N^t(W_{m+t}) = W_{m-t}$.
\\ \\
3) As a biproduct of 2) it follows that $W_{m+t-1}/W_{m-t-1}\subset \Ker(N^t|_{W_{m+t}/W_{m-t-1}})$. For the other inclusion it is enough to prove  that $(N^t)^{-1}(W_{m-t-1})\cap W_{m+t} \subset W_{m+t-1}$.
\\ \\
Indeed, since $W_{m-t-1}=\sum_{a} M_{m+t-a,a}$ with $0\le m+t-a\le m$,  it is enough to prove that $(N^t)^{-1}(M_{m+t-a,a})\cap W_{m+t} \subset W_{m+t-1}$. 
Observe that $0\le m+t-a\le m$ if and only if $a-m\le t\le a$. \\ \\
If $z\in (N^t)^{-1}(M_{m+t-a,a})\cap W_{m+t}$, then 
\[N^t(z)\in \Ker N^{m-(m+t-a)}\cap \im N^a = \Ker N^{a-t}\cap \im N^a.\] Hence,  $N^{t}(z) = N^{a}(w) $  or  $N^{t} (z - N^{a-t}(w)) = 0$ for some $w$, i.e., \[z\in \Ker N^t + \im N^{a-t}.\]
But  $N^t (z) \in \Ker N^{a-t} $, then  $ 0 = N^{a-t}( N^t (z)) =  N^a (z) $, i.e.  $ z \in \Ker N^a $.
Therefore \[ z  \in ( \Ker N^t + \im N^{a-t} ) \cap \Ker N^a = \Ker N^t + \Ker N^a\cap \im N^{a-t}, \]  hence  $ z  \in M_{m-t,0} + M_{m-a, a-t}\subset W_{m+t-1}$.
\hfill \qed
\ \\ \ \\
In order to compute the limit Hodge structure for the VPHS $\widetilde{\sH}^3$ at a point $p\in\PP^1\backslash B$, we apply the formula obtained in example \ref{app} and the fact that  $ N^r : {\Gr}_{n+r}^W   \simeq  \Gr_{n-r}^W $ for all $r$ for the weight filtration centered at $n$ (we say it is {\it symmetric} at n ). Let us define $ n_i \stackrel{def}{=} \dim ( \im N_i)$ and $ m_i \stackrel{def}{=} \text{dim} (\Ker N_i) $ for  $ i \in \{ 1,2,3 \} $. To simplify the notation of the proof  of the following proposition we omit the tildes in  the  components $H^{i,j}$ of the local system $\widetilde{\sH}^3$.
\begin{prop}
\label{lmhs}
The limit Hodge filtration $(W_{\ldotp},F^{\cdotp}_{\infty})$  for the family $\tilde{\sX}$ can be described as follows, where $N_{i,j} \stackrel{def}{=} N_i \circ N_j $ and 
$ o \stackrel{def}{=}  \dim ( M_{2,1} )$:
\begin{enumerate}
 \item  $N= 0 $ and it is pure of weight three.
 \item $N \neq 0, N^2 = 0 $ there are two cases:
       \begin{enumerate} 
        \item $N_1 \neq 0, N_3 \neq 0$ such that $ N_{2,1} = N_{3,2} = 0 $,
         \item $N_1=N_3= 0 $ with $N_2 \neq 0 $.
         \end{enumerate}
         For these cases  the weight filtration centered at three is:
          \[
\Gr_i^W( H_{\Q} ) = \left\{ \begin{array}{crr}

                                        0     \hfill         &    i=0,1,5,6   \\
                                         \C^{ n_2 + 2}  \hfill     & i = 2,4 & \mbox{a)},  \\
                                         \C^{n_2}     \hfill             &  i=2,4  & \mbox{b)}, \\
                                        \C^{2(m_2-1)}  \hfill         &  i=3  &\text{a)}, \\
                                        \C^{2(m_2 +1)} \hfill       & i= 3  & \mbox{b)}. \\
                                         \end{array}
                                        \right.
\]
 \item $N^2\neq 0, N^3= 0 $.  
 \[
\Gr_i^W( H_{\Q} ) = \left\{ \begin{array}{cr}
                                        0  \hfill      & i = 0,6  \\
                                         \C^2        \hfill          &  i=1,5 \\
                                        \C^{o -2}  \hfill         &  i=2,4 \\
                                        \C^{2(k+1-o)} \hfill       & i= 3.  \\    
                                        \end{array}
                                        \right.
\]
 \end{enumerate}
\end{prop}
\pf 
\begin{enumerate}
 \item $N= 0 $ ; the weight filtration centered at three is: $ W_i = 0 $ for  $ i \in \{ 0,1,2 \} $ and otherwise $ W_j = H_{\Q} $.
 Therefore, $ \Gr_3 (H_{\Q} ) = H_{\Q} $. 
 \item  Assume that $N \neq 0, N^2 = 0 $. We have the following  general decompositions for $ \{ \Ker N^i, \im N^i  \}_{ i=1,2} $:
\[ 
\begin{matrix}
 \Ker N = \oplus_{i=1}^3 \Ker N_i \oplus H^{0,3} \hskip0.3cm{,}  \hskip1.9cm{} \im N =  \oplus_{i=1}^3 \im N_i , \cr 
 \Ker N^2 =  \Ker N_{2,1} \oplus  \Ker N_{3,2}   \oplus H^{1,2} \oplus H^{0,3}, \hskip0.7cm{}  \im N^2 = \im N_{2,1} \oplus \im N_{3,2}\hfill  
 \end{matrix}
\]
and in both cases the weight filtration is
 given as:
\begin{eqnarray*} 
 W_0 = W_1 = 0, W_2 = \im N, W_3 = \Ker N , W_4 = \Ker N^2 = W_5 = W_6 = H_{\Q}. 
 \end{eqnarray*}  
 In particular we get $Gr^W_iH_{\Q}=0$ for $i\in\{0,1,5,6\}$ in both cases, as stated.
 \\ \\
For the other graded groups we have:
\\ \\
 \noindent a) $N_1 \neq 0$ (and hence $N_3 \neq 0$ since  the polarization is non-degenerate and the GM-connection is 
 compatible with the metric induced by it)
 with 
$ N_{2,1} = N_{3,2} = 0 $
 . Hence  
 the weight filtration simplifies further to:
 $$ 
 0 = W_0 = W_1 \subset \im N \subset \Ker N  \subset W_4 = W_5 = W_6 = H_{\Q} 
 $$ 
 
  \noindent b) $N_1 = N_3 = 0 $ with $N_2  \neq 0 $. In this case:
  \begin{eqnarray*}
\im N_2 = \im N \subset \Ker N = H^{3,0} \oplus H^{0,3} \oplus H^{1,2} \oplus \Ker N_2 .
\end{eqnarray*}
Then for both cases above:

\[
W_3  = \Ker N  = \left \{ \begin{array}{cr}
                          \Ker N_2 \oplus \Ker N_3  \oplus  H^{0,3}    &           \text{ a)}, \\
                          H^{3,0}  \oplus H^{0,3}  \oplus  H^{1,2} \oplus \Ker N_2&  \text{ b)}. \\
                          \end{array}
                                        \right.
\]

\[
 W_2 =   \Gr_2^W H_{\Q} =\im N  = \left \{   \begin{array}{cr}
                       \oplus_{i=1}^3 \im N_i = \C  \oplus  \C^{n_2} \oplus  \C      &  \text{  a) }, \\
                       \im N_2   = \C^{n_2}                  & \text{ b) } .\\
                      \end{array}
                         \right.
\]
Trivially:
\begin{eqnarray*}
 \Ker N_2 / \im N_1  \simeq \C^{m_2-1} \simeq \Ker N_3 / \im N_2 ,  \, \, H^{1,2}  / \im N_2 \simeq \C^{m_2}.
\end{eqnarray*}
Hence:
\[
 \Gr_3^W H_{\Q} = \left \{   \begin{array}{cr}
                       \Ker N_2 / \im N_1 \oplus  \Ker N_3 / \im N_2 \simeq      \C^{2(m_2-1)}  &  \text{  a) }, \\
                         H^{3,0}  \oplus H^{0,3}  \oplus  H^{1,2}/ \im N_2 \oplus \Ker N_2 \simeq  \C^{2(m_2+1) } & \text{ b) }. \\
                      \end{array}
                         \right.
\]
\ \\ 
\item  $ N^2 \neq 0,  N^3 = 0 $ with $ N_{2,1}  \neq 0, N_{3,2}  \neq 0 $. 
The weight filtration is explicitely:
$$
\begin{array}{crl}
W_0 & =  &   0, \\
W_1 &  =   &   \im N_{2,1} + H^{0,3} \simeq \C^2 , \\
W_2 & =  &  \im N \cap \Ker N  \simeq \C^o \, (\text{note}: M_{1,2} \subset M_{2,1}) ,\\
W_3  & = & \im N + \Ker N  \simeq \C^{2k+2 -o}, \\
W_4 & = & \Ker N^2  = \Ker N_{3,2} + H^{1,2} + H^{3,0} = \C^{k-1} \oplus \C^k \oplus \C^1 = \C^{2k}, \\
W_5 & =  &  W_6 = H_{\Q}.
\end{array}
$$
from which:
\[
\Gr_k^W H_{\Q}   =    \left \{   
                        \begin{array}{cr}
                        0    &  \text{for } \, k=0 , \\
                       \C^2  &  \text{for } \, k = 1 ,\\
                         \C^{o-2} & \text{for }  \, k = 2 , \\
                         \C^{2(k+1-o)}      &  \text{for }  \, k = 3. \\
                      \end{array}
                         \right.
\]
\end{enumerate}

$\text{ }$\hfill \qed

%%%%%%%%%%%%%%%%%%%%%%%%%%%%%%%%%%%%%%%%%%%%%%%%%%%%%%%end%%%%%%%%%%%%%lmhs%%%%%%%%%%%%%%%%%%%%%%%%%%%%%%%%%%%%%%%%%%%%%%%%%%%%%%%%%%
\ \\
\noindent  We return to the study of the VMHS for the family $ \sX$  over $B$ considered in the introduction.
 \\ \\
 By assumption, the family $ \sX \stackrel {{f}}\longrightarrow B\subset\PP^1$ is real analitically 
 trivial, i.e. the sheaf $R^3f_*\C$ is a local system on $B$. Additionally, by  the RH-correspondence, there exist a GM
 -connection $\xymatrix{\sH^3\ar[r]^(0.4){\nabla^{GM}} & \sH^3\otimes\Dif_B^1}$. Moreover, the weight filtration on the fibers 
 fits together to form a subbundle $\sW_2R^3f_*\Q\subset R^3f_*\Q$ and we have a short exact sequence
 (see also equation (\ref{Clemens-Schmidt})):
 \[
 \xymatrix{
 0\ar[r] & \sW_2R^3f_*\Q\ar[r] & R^3f_*\Q\ar[r]^{{\pi}^*} & R^3\widetilde{f}_*\Q\ar[r] & 0
 }
 \]
 \noindent
 A trivialization for $R^3f_*\C$ induces a trivialization for $\sW_2R^3f_*\C$ and so the action of the monodromy on $R^3f_*\C$ is
 compatible 
 with the action of the monodromy on $\sW_2R^3f_*\C$, in particular the GM-connection 
 on $\sW_2R^3f_*\C\otimes\sO_B=\sW_2\sH^3$ is just 
 the restriction of $\nabla^{GM}$ on $\sH^3$ to $\sW_2\sH^3$ and by passing to the quotient
 the short exact sequence above induces a 
 connection $\overline{\nabla}^{GM}$ on $\widetilde{\sH}^3$ with flat sections $R^3\widetilde{f}_*\C$. By the uniqueness of
 the  GM-connection (see proposition 2.16 on \cite{EDPSR}), this connection 
 is none other than  $\widetilde{\nabla}$ on $\widetilde{\sH}^3$, i.e. we have a short exact sequence which is 
 compatible with the GM-connection:
 \begin{eqnarray}
 \label{compatibility}
 \xymatrix{
 0\ar[r] & \sW_2\sH^3\ar[d]^{\nabla^{GM}}\ar[r] & \sH^3 \ar[d]^{\nabla^{GM}}\ar[rr]^{\widetilde{\pi}^*} & & \widetilde{\sH}^3 \ar[d]^{\widetilde{\nabla}}\ar[r] & 0 \cr
 0\ar[r] & \sW_2\sH^3\otimes\Dif_B^1 \ar[r] & \sH^3\otimes\Dif_B^1\ar[rr]^{\widetilde{\pi}^*\otimes id} & &\widetilde{\sH}^3\otimes\Dif_B^1\ar[r] & 0
 }
 \end{eqnarray}
 \begin{prop}
 \label{locsys}
 $(\sH^3,\nabla^{GM})$ is a VMHS.
 \end{prop}
 
 \pf Observe that the Hodge filtrations are compatible, so we have a commutative diagram with exact arrows
 \[
 \xymatrix{
 0\ar[r] & 0\ar@{=}[d]\ar[r] & \sF^3\sH^3 \ar@{^{(}->}[d]\ar[r]^{\widetilde{\pi}^*} &  \sF^3\widetilde{\sH}^3\ar@{^{(}->}[d]\ar[r] & 0 \cr
 0\ar[r] & 0\ar@{^{(}->}[d]\ar[r] & \sF^2\sH^3 \ar@{^{(}->}[d]\ar[r]^{\widetilde{\pi}^*} &  \sF^2\widetilde{\sH}^3\ar@{^{(}->}[d]\ar[r] & 0 \cr
 0\ar[r] & \sF^1\sW_2\sH^3\ar@{=}[d]\ar[r] & \sF^1\sH^3 \ar@{^{(}->}[d]\ar[r]^{\widetilde{\pi}^*} &  \sF^1\widetilde{\sH}^3\ar@{^{(}->}[d]\ar[r] & 0 \cr
 0\ar[r] & \sW_2\sH^3\ar[r] & \sH^3 \ar[r]^{\widetilde{\pi}^*} &  \widetilde{\sH}^3\ar[r] & 0
 }
 \]
 where all the rows are exact. 
 In particular $\nabla^{GM}(\sF^p\sH^3)\subset\sF^{p-1}\sH^3\otimes\Dif_B^1$  by the commutativity of diagram (\ref{compatibility}), hence it becomes a VMHS. \hfill \qed
 \\ \\
 \noindent Since $ f$ is quasi-projective this VMHS is in fact graded-polarizable, indeed this is a {\it geometric variation of  mixed Hodge structure.}
 We claim even that
 \begin{cor}
  $f$ is  a geometric VMHS and an admissible  variation of Hodge structure in the sense of Steenbrink-Zucker ( see  \cite{peters} Theorem 14.51  and  \cite{stzu} ).
 \end{cor}

%%%%%%%%%%%%%%%%%%%%%%%%%%%%%%%%%%%%%%%%%%%%%%%%%%%%%%version diez%%%%%%%%%%%%%%%%%%%%%%%%%%%%%%%%%%%%%%%%%%%%%% 
A desingularization of the family produces a VPHS  $\widetilde{\sH}^3$ whose limit MHS can be described as in proposition \ref{lmhs}. 
\ \\ \\
 Denote by $\M_l ( \C)$ the set of $l$ by $l$ matrices over $\C$ and denote by $J(l) \in  \M_l ( \C)$  
 the  Jordan matrix with entries
  \[ 
            J(l)_{r,s}  =  \left \{   
                        \begin{array}{cr}
                        1   &  \text{for } \,  r = s + 1, \\
                         0      &  \text{otherwise}.
                      \end{array}
                         \right.  
  \]
 Then : $\rank \Ker (J(l)^i ) = i $ for $ i \le l $, $\rank \Ker (J(l)^{i+1}) - \rank \Ker (J(l)^i)  = 1 \;\forall i<l$.
The Jordan form ${\bf J} (A)$ of a nilpotent matrix $ A \in \mbox{M}_{n}( \C) $  
is written as a direct sum of the corresponding
Jordan  block matrices. We call such a direct sum of 
Jordan block matrices simply a {\it Jordan matrix } of a Jordan form.  If a Jordan block matrix $J(m)$ appears with 
multiplicity $r$ we denote it by 
$J(m)^r$.

Assume we have a VHS of type $(1,k,k,1)$. Keeping the notation above for $T$, $T_u$ and $N$,  if $N\ne 0$ and $k\ge 2$,  then the  Jordan matrix of $N$ is one 
of the following types:
\begin{itemize}
 \item  Type (1): $J(4) \oplus J(2)^s $ with  $ s \leq k-1$,
 \item  Type (2): $ J(3)^2 \oplus J(2)^s $ with $ s \leq k-2$,
 \item  Type (3): $J(2)^s $ where $ s \leq k +1 $,
 \item  Type (4): $ J(3)  \oplus J(2)^s $ with $ s \leq [\frac{ 2k-1}{2}]$.
\end{itemize}

\begin{prop}
\label{jdec}
The Jordan canonical form of $N$ is of type (1),(2) or (3).
 \end{prop}

\pf
\begin{enumerate}
\item A type  (3)  Jordan matrix decomposition implies that there are at most $ k + 1$ two by two blocks.
This implies that $N^2 = 0 $.
\item A type (1) Jordan matrix decomposition correspond to the maximal unipotent case, 
which is known to occur for instance for the family of \cite{candelas}.
\item If $N^3=0$ but $N^2\ne 0$, we know from linear algebra that all Jordan blocks are of size 3, 2 or 1, 
which correspond to either type (2) or type (4). \\
  A type (4)  Jordan Matrix decomposition is not possible. 
 For that recall the abstract situation of example~\ref{jord},~namely:
\begin{lem}
\label{corta} Recall the notation of  proposition \ref{lmhs}:
 \begin{eqnarray*}
  H^{3,0}  \stackrel{N_1}{\rightarrow} H^{2,1} \stackrel{N_2}{\rightarrow} H^{1,2} \stackrel{N_3}{\rightarrow} H^{0,3}
  \end{eqnarray*}
where $ H^{3,0} \simeq  H^{0,3} \simeq \C, H^{2,1} \simeq H^{1,2} \simeq  \C^{k} $ then:
 $ N_1$ is one-to-one  $\Leftrightarrow$  $ N_3 $ is surjective.
 \end{lem}
 \pf  (of the lemma). The polarization $Q$ is flat with respect to the connection  $N $. 
 \hfill \qed \ \\
 \item
Assume $N^3\equiv0$ but there exist a  three-dimensional $N$-cyclic space  
 \[ 
 W_0=<w,N(w),N^2(w)>.
 \] 
 Without loss of generality, assume either $w\in H^{3,0}$ or $w\in H^{2,1}$. Indeed, write $w= v_0+v_1+v_2+v_3$ with $v_j\in H^{k-j,j}$. Then
 \[
 N(w) = N(v_0) + N(v_1) + N(v_2)
 \]
 and
 \[
 N^2(w) = N^2(v_0) + N^2(v_1), 
 \]
 since $N(v_3)=N^2(v_2)=0$.
 \ \\ \ \\
If $N^2(v_0)\ne 0$, then $<v_0, N(v_0), N^2(v_0)>$ is a three-dimensional N-cyclic space. On the other hand, if $N^2(v_1)\ne 0$, then
$<v_1, N(v_1), N^2(v_1)>$ is a three-dimensional N-cyclic space.
\\ 
 \begin{enumerate}
 \item If $w=v_0\in H^{3,0}$, since $Q$ is non-degenerate there exist a $u\in H^{2,1}\backslash\{0\}$ 
 such that $Q(u,N^2(w))=1$ and because of $Q$-flatness of the VHS  with respect to $\nabla^{GM}$:
 \[
 Q(N(u), N(w)) + Q(u, N^2(w))=0
 \]
 therefore $ Q(N(u),N(w))=-Q(u,N^2(w)) = -1$  thus $N(u)\ne 0$. 
 \\
 Similarly, 
 $Q(N^2(u), w)= - Q(N(u),N(w)) = 1$ and $N^2(u)\ne 0$ either, therefore
  $W_1:=<u,N(u),N^2(u)>$ is a another three-dimensional $N$-cyclic space and
 $W_0\cap W_1=0$.
 \\
 \item  If $w=v_1\in H^{2,1}$ then $N^2(w)\in H^{0,3}\backslash\{0\}$ and  because of the non-singularity of $Q$, there exist a $u\in H^{3,0}\backslash\{0\}$
such that $Q(u,N^2{w})=1$. As before, we will have 
\[
 Q(N(u), N(w)) + Q(u, N^2(w))=0
 \]
 therefore $ Q(N(u),N(w))=-Q(u,N^2(w)) = -1$  thus $N(u)\ne 0$ and again, $Q(N^2(u), w)= - Q(N(u),N(w)) = 1$ and $N^2(u)\ne 0$ either, therefore
  $W_1:=<u,N(u),N^2(u)>$ is a another three-dimensional $N$-cyclic space and
 $W_0\cap W_1=0$.
 \hfill \qed
 \end{enumerate}
\end{enumerate}

\bibliographystyle{plain}
\renewcommand\refname{References}

{}

\end{document}